\documentclass[11pt]{article}

\usepackage{amsmath,bbm,amsfonts,latexsym,amssymb}
\usepackage{graphicx}
\usepackage{enumerate}
\usepackage{natbib}
\usepackage{url} 
\usepackage{fancyhdr}
\usepackage{xcolor}

\usepackage{algorithm}
\usepackage{algorithmic}
\usepackage{hyperref}

\newtheorem{theorem}{Theorem}[section]
\newtheorem{lemma}[theorem]{Lemma}
\newtheorem{proposition}[theorem]{Proposition}
\newtheorem{corollary}[theorem]{Corollary}

\newtheorem{definition}[theorem]{Definition}
\newtheorem{remark}[theorem]{Remark}

\newtheorem{example}[theorem]{Example}

\newcommand{\calB}{{\cal B}}
\newcommand{\calC}{{\cal C}}

\newcommand{\calH}{{\cal H}}

\newcommand{\calN}{{\cal N}}

\newcommand{\calS}{{\cal S}}

\newcommand{\veps}{{\varepsilon}}
\newcommand{\eps}{{\epsilon}}

\newcommand{\la}{\lambda}

\def\E{{\mathbb{E}}}

\def\eps{\epsilon}

\newcommand{\real}{{\mathbb{R}}}
\newcommand{\pr}{{\mathbb{P}}}
\newcommand{\qed}{$\Box$}

\graphicspath{{hd-plots/}}

\date{}

\begin{document}

    \title{
   Discriminating Tail Behavior Using Halfspace Depths: Population and Empirical Perspectives   
        }
        \author{Sibsankar SINGHA\textsuperscript{(a)} \\[.5ex]
        Marie KRATZ\textsuperscript{(b)} and Sreekar VADLAMANI\textsuperscript{(c)}
        \\[1ex]
        \small
        \textsuperscript{(a)} T\'el\'ecom Paris, Palaiseau, France; Email: sibsankar.singha@telecom-paris.fr \\
        \small
        \textsuperscript{(b)} ESSEC Business School, CREAR, Cergy-Pontoise, France; Email: kratz@essec.edu \\
        \small
        \textsuperscript{(c)} TIFR-CAM, Bangalore, India; Email: sreekar@tifrbng.res.in
        }
        
\maketitle

\begin{abstract}
   We study the empirical version of halfspace depths with the objective of establishing a connection between the rates of convergence and the tail behaviour of the corresponding underlying distributions. The intricate interplay between the sample size and the parameter driving the tail behaviour forms one of the main results of this analysis. The chosen approach is mainly based on weighted empirical processes indexed by sets by Alexander (1987), which leads to relatively direct and elegant proofs, regardless of the nature of the tail. 
This method is further enriched by our findings on the population version, which also enable us to distinguish between light and heavy tails. These results lay the foundation for our subsequent analysis of the empirical versions. 
Building on these theoretical insights, we propose a methodology to assess the tail behaviour of the underlying multivariate distribution of a sample, which we illustrate on simulated data. 
The study concludes with an application to a real-world dataset.
\end{abstract}


\newpage 

\section{Introduction}
 
Recent progress in high-dimensional statistics has led to a renewed focus on developing tools that clarify the geometric structure of datasets. A variety of multivariate quantiles and statistical depth functions have been introduced, offering nonparametric insights into data in multiple dimensions. These tools have proven to be particularly valuable for addressing statistical inference challenges. 
 
In contrast to quantiles, which are determined analytically through the inverse of the cumulative distribution function, depth functions embrace a geometric perspective; they utilize halfspaces, paraboloids, and projections to assess centrality from a wider viewpoint, resulting in an arrangement of observations that extends outward from the center.

Numerous depth functions have been introduced and studied, starting with Mahalanobis distance depths (\cite{Mahalanobis2018,Liu1993,Zuo2000}), the well-known and used Tukey or halfspace depth (\cite{Tukey1975}), going on, for instance, with simplicial (volume) depths (\cite{Oja1983,Liu1990}), onion depths (\cite{Barnett1976,Eddy1982}), all notions of spatial depths (\cite{Dudley1992,Chaudhury1996,Koltchinskii1997,Vardi2000,Mottonen2005}), the projection depth (\cite{Donoho1992,Zuo2003,Dutta2012,Nagy2020}), the zonoid depth (\cite{Dyckerhoff1996,Koshevoy1997,Koshevoy2003}), 
local depths (\cite{Agostinelli2011,Paindaveine2013}). We refer to \cite{Hallin2010,Mosler2002,Mosler2013,Kuelbs2016,Chernozhukov2017,Nagy2020,Nagy2021,Mosler2021} for theoretical and practical aspects (as well as computational) of depth functions, as well as to \cite{Nagy2022} and references therein for halfspace depths, and \cite{Nagy2025} for a more recent analysis. 
Our focus is on halfspace depth. 
Many results are already available for the population-based analysis of this measure, but less for its empirical estimation, particularly in terms of its asymptotic (extreme) behaviour; see e.g. the work by \cite{He2017} for multivariate regularly varying distributions.
Considering practical applications, the questions regarding asymptotics become even more critical when examining sample versions of this measure. We aim to identify tail properties of the underlying probability measure from the asymptotic behaviour of the halfspace depth. 

To do so, we first investigate the extreme behaviour of halfspace depth based on the nature (light or heavy tail) of the underlying distribution, establishing decay rates for, both, population and sample versions.
 
Recall that (see \cite{Donoho1992}), as the sample size increases, the halfspace depth for a sample converges almost surely to the halfspace depth for the underlying distribution. To obtain rates of decay of halfspace depth,  regardless of the nature of the tail, light or heavy, our approach builds mainly on the theory of empirical processes (\cite{Shorack2009}) and weighted empirical processes indexed by sets by \cite{Alexander1987}. The latter paper and its powerful results, helped us prove our results in a rather direct and elegant way. 

For heavy tails, we highlight related works by \cite{Einmahl2015} and \cite{He2017}, who propose an estimator for half-space depth that extends beyond the data hull using extreme value theory.  
Two different approaches are developed within the framework of (multivariate) regular variation ( (M)RV), the definition of which is provided in Appendix~\ref{app:He2017}. In \cite{Einmahl2015}, the authors enhance the estimation of empirical half-space depth in the tail by extending it beyond the data hull, using the limiting extreme value distribution for the one-dimensional tail probability along each projection direction. In \cite{He2017}, a direct multivariate extreme value approach is considered to estimate extreme quantile regions.

In the scenario where the tail of the distribution is light, \cite{Burr2017} employed a novel geometric methodology to derive uniform convergence rates for halfspace depth. Through the reorganization of halfspaces into one-dimensional family, the authors successfully improved convergence bounds for the sample version of halfspace depth, surpassing the typical Glivenko--Cantelli bounds. This advancement was particularly evident when considering exponential decay in the underlying distribution.

Our approach enables the establishment of convergence results for both light and heavy tails, within the data hull. It further investigates whether depth functions reflect the tail behaviour of the underlying distribution, and examines the role of sample size in visualizing this behaviour. We then use the obtained results to address the inverse problem of identifying the tail characteristics of the underlying distribution through the asymptotics of the empirical halfspace depth of a given dataset.

Our main results are illustrated on simulated data. We apply the developed methodology on a real--world dataset.
The computation of halfspace depth is performed using the data-depth Python library. For depth contours, we use the R package ddalpha developed in \cite{Pokotylo2019}, and also refer to the R package TukeyRegion developed in \cite{Liu2019}. For faster computation on large datasets, we employ an approximation of the true depth, as proposed by \cite{Dyckerhoff2004} (see also \cite{Dyckerhoff2021}), which is also available in both the data-depth and ddalpha packages. Convergence rates for this approximation have been established by \cite{Nagy2020}. We note that the computation of depth functions remains a significant challenge, actively addressed by various research teams who continue to develop efficient and scalable software (see, e.g., \cite{Genest2019, Liu2015, Mahalanobish2015, Fojtik2024}).

 {\it Notation.} All the analysis in this paper is on $\real^d$ equipped with the Borel sigma algebra $\calB(\real^d)$, unless otherwise stated. The centered unit open ball and the unit sphere in $\real^d$ are denoted by $B^d$ and $S^{d-1}$, while $\langle\cdot,\cdot\rangle$ and $\|\cdot\|$, denote the Euclidean inner product and $\ell^2$-norm, respectively, in $\real^d$.

{\it Structure of the paper.} 
Section~\ref{sec:multi-quantiles-HD} considers the population version of halfspace depth, completing the literature on the topic with insightful results for the study of its empirical counterpart. The latter is developed in Section~\ref{sec:empirical-HD}. In particular, we investigate the asymptotic behaviour of the sample version in relation to the sample size. These results form the foundation for developing a methodology in Section~\ref{sec:HDmethod} to identify the nature of the tail behaviour of the underlying measure, which is then applied to a real--world dataset. 
All the proofs are presented in Section~\ref{sec:proofs}, with necessary supplementary material provided in Appendices~\ref{app:He2017} and \ref{subsec:app2}. 
Additional illustrative examples are given in Appendix~\ref{app:illustration}.


\section{Halfspace depth} 
\label{sec:multi-quantiles-HD}

As in the univariate case, it is expected that depth functions encode the tail behaviour of the underlying probability measure. It is this line of thought that we explore in this section, studying the asymptotic behaviour of halfspace depth {\it vis-\`{a}-vis} of the tail characteristics of the underlying probability measure, while also discussing the {\it status quo} of the subject to contextualise our results.

Recall that a {\it depth function} corresponding to a probability measure $\pr$ is a non-negative function $D(x,\pr)$ defined at every point $x \in \real^d$, which provides an outward ordering from the {\it center of the distribution}.

Since it is desirable that depth functions decrease to zero in every direction from the median/center, the {\it center} of a distribution is defined as the point of maxima of the depth function. In case the maximum is attained at multiple (finitely many) points then the centroid of all such points is called the median (and the center). 

We shall focus on the following notion of halfspace depth. It is widely used, and is also a good representative of the class of depth functions, as it satisfies most of the desirable properties for depth functions (see \cite{Mosler2021}, Table 2). 
\begin{definition}[Halfspace depth; \textnormal{\cite{Tukey1975}}] 
    \label{def:HD}
    For a probability distribution $\pr$ defined on $\real^d$, the halfspace depth is given by: 
    $HD(x,\pr)= \inf \{\pr(H):\,H \in \calH_x\},$ where $\calH_x$ denotes the set of halfspaces in $\real^d$ containing $x\in\real^d$. Specifically, if $\pr$ has a probability density function $f$, then
$\displaystyle
HD(x,\pr)=\inf_{\|p\|=1}\int_{\{y: \langle y-x,p\rangle \geq 0\}} f(y) \,  dy.
$
\end{definition}
Intuitively, a high depth point indicates that it is more central, while a low depth point denotes a relatively extremal point.

\begin{figure}[h]
    \centering
     \includegraphics[trim=0cm 0cm 0cm 2cm, clip, width=.6\linewidth]{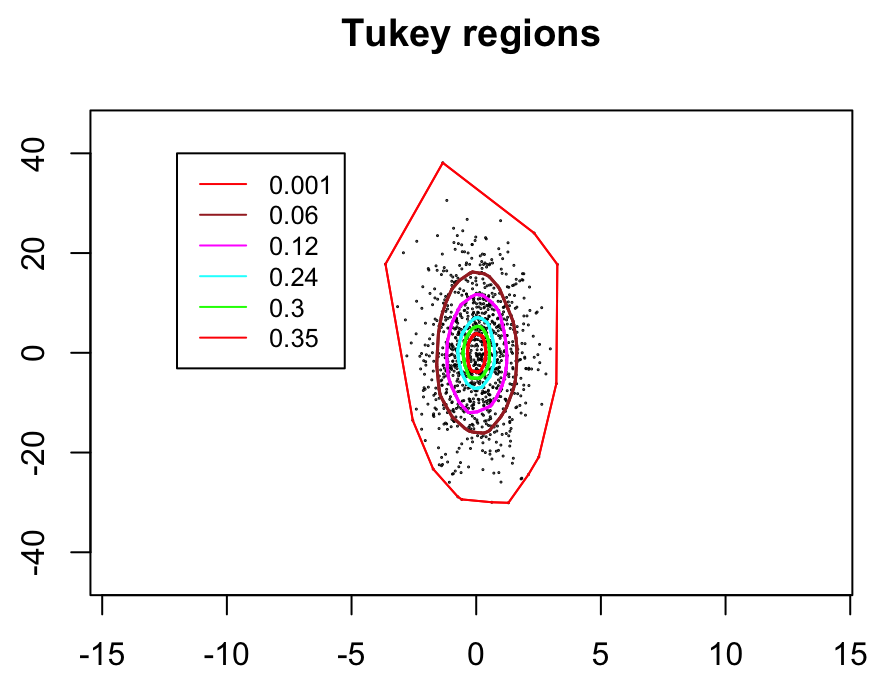}
     \vspace{-2.5ex}
     \parbox{350pt}{\caption{\sf Representation of the Tukey depth contours for $6$ different depths, considering a sample of $1000$ observations (black points) from a mean zero Gaussian distribution with covariance $\text{diag}(1,100)$.}
     \label{fig:Tukey-regions}}
 \end{figure}

This is illustrated in Figure~\ref{fig:Tukey-regions} when considering a Gaussian sample (using the R-Package from \cite{Barber2022}, based on \cite{Liu2019}). Note that the symmetry present in the underlying sampling distribution may not be seen precisely in the simulated sample. Isodepth contours are drawn for different values of depth. Observe that the asymmetry is more evident in the extremes, than in the bulk region.

\vspace{1ex}


Continuing on the theme of identifying symmetry, we state the following result which establishes that the halfspace depth corresponding to a probability measure  does inherit (asymptotic) symmetry of the underlying measure, whenever the underlying parent measure has density with respect to Lebesgue measure.

\begin{theorem}\label{thm:asymp-depth-ell-symm}
If the probability density function $f$ of a distribution $\pr$ satisfies $\displaystyle \lim_{t\to\infty} \sup_{\|\Sigma x\|=\|\Sigma y\|} \frac{f(tx)}{f(ty)} = 1$, where $\Sigma$ is a symmetric and positive definite matrix, then the halfspace depth for $\pr$ satisfies 
$$ 
\lim_{t\to\infty} \sup_{\|\Sigma x\|=\|\Sigma y\|}\frac{HD(tx,\pr)}{HD(ty,\pr)} =1.
$$
\end{theorem}

{\it Decay rate of halfspace depth:} 
Recall that our motivation is to understand the connection between the extremal behaviour of a probability measure and the asymptotics of halfspace depth functions. 
Specifically, we question the halfspace asymptotic decay depending on the light or heaviness of the underlying probability measure. We shall start with a generic result comparing distributional asymptotics with that of the induced halfspace depth, when marginals of the joint distribution may have different asymptotic behaviours. For that, we first recall another equivalent way to define the halfspace depth of a distribution in terms of its marginals (see e.g. \cite{Donoho1992}), as follows:
\begin{equation}\label{eqn:donoho-depth}
    HD(x,\pr) = \min_{h:\|h\|=1}(1- F_h(\langle h,x\rangle))
\end{equation}
where $F_h$ is the c.d.f. of the univariate projection of $\pr$ onto the $h$ direction. In order to explicitly highlight the roles of left and right tails, we rewrite relation \eqref{eqn:donoho-depth} using the relative positioning of $x$ and $h$. Specifically,
\begin{equation}\label{eqn:hdepth-project}
   HD(x,\pr) = \min_{\substack{ \|h\| =1 \\ \langle h,x\rangle \ge 0}} \Big(\min \Big\{(1- F_h(\langle h,x\rangle)),\, F_h(\langle h,x\rangle)\Big\} \Big),
\end{equation}

which shall be the core ingredient of our analysis throughout the paper.
As a direct consequence of \eqref{eqn:hdepth-project}, we can state the following simple but interesting observation, which was noticed by \cite{Dyckerhoff2004}. This will turn out to be very useful for deriving a method to discriminate between light and heavy tails (see Section~\ref{sec:HDmethod}).

\begin{proposition}
\label{thm:hd-decay-survivalMargins}
Let $\pr$ be a probability measure on $(\real^d,\calB(\real^d))$, and let $\{e_i\}_{i=1}^d$ be any orthonormal basis of $\real^d$, $F_{e_i}$ denoting cumulative distribution function of the marginal of $\pr$ along $e_i$ direction.
Then, we have the following upper bound:
\begin{equation}
\label{eq:gl-upperBd-HD}
HD(t\,x,\pr)\leq \min\left(\min_{1\le i\le d} \left(1-F_{e_i}(t \, \langle x,e_i\rangle)\right),\min_{1\le i\le d} F_{e_i}(t \, \langle x,e_i\rangle)\right).
\end{equation}
In particular, by choosing $\{e_i\}_{i=1}^d$ as the canonical basis, the upper bound in \eqref{eq:gl-upperBd-HD} corresponds to the minimum over all the marginals $F_i$ along the canonical basis:
$$HD(t\,x,\pr)\leq \min\left(\min_{1\le i\le d} \left(1-F_{i}(t\,x_i)\right), \min_{1\le i\le d} F_{i}(t \, x_i) \right),
$$
where $x_i = \langle x,e_i\rangle$ denotes the projection of $x$ on $e_i$.
\end{proposition}
As a consequence, if the marginal distribution along $e_i$ has a right (left) light tail for some $1\le i\le d$, then the decay of halfspace depth cannot be slower than exponential along $x$ directions satisfying $\langle x,e_i \rangle >0$ ($<0$). If all marginals are heavy-tailed, then the bound will correspond to the one with the smallest tail index. This offers a simple tool to discriminate between light and heavy tails, which we will develop further in Section~\ref{sec:HDmethod}. 

Next, we investigate a lower bound asymptotic decay of halfspace asymptotic depth, specifically considering the tail behaviour of the underlying probability measure:
\begin{theorem} 
\label{thm:hd-decay-tail}
    Assume there exists a positive function $g$ such that 
    \begin{equation}\label{eqn:unif-direction-density-decay}
        \liminf_{R \to \infty} \inf_{\theta\in S^{d-1}}\dfrac{\pr\big( \langle\theta, X\rangle \geq R\big)}{g(R)} >0.
    \end{equation}
Then, for $x\neq 0$, there exists $t_0>0$ such that, 
    \begin{equation*}
        HD(t x, \pr) \geq c\,\, g(t\|x\|), \,\,\,\,\, \forall \,\,t > t_0,
    \end{equation*}
    for some constant $c>0$ independent of $x$ and $t$.
    
As an immediate consequence, we have:
    $$
    \pr(H) \ge c\, g(t\|x\|) , \qquad\forall \,\,H\in\calH_{tx}.
    $$
\end{theorem}

Now, we question the halfspace asymptotic decay depending specifically on the light or heaviness of the underlying probability measure. The heavy tail case has been studied in \cite{He2017}, providing the decay rate in \cite[Proposition~2]{He2017} for multivariate regularly varying distributions (under some additional conditions). 
We recall the latter in Appendix~\ref{app:He2017},  
also adding another version under the condition that the probability measure has a density.
Now, we aim at completing the picture by providing in Theorem~\ref{thm:hd-decay-light-tail} the corresponding result for the light tail case.
\begin{theorem} 
\label{thm:hd-decay-light-tail}
Let $\pr$ be a probability measure on $(\real^d,\calB(\real^d))$.
\begin{itemize}
\item[(i)] Let the moment generating function $M_{\pr}(h) = \int_{\real^d}e^{\langle h,y\rangle}\,\pr(dy)$ be finite for all $h$ in an open neighbourhood $\calN$ of $0$. Then, the halfspace depth is also light tailed, i.e., for $x \neq 0$, 
\begin{equation}\label{eqn:light-tail}
    \liminf_{t\to\infty} \frac1t \log \left( \frac1{HD(tx,\pr)}\right) >  \sup_{h\in\calN} \langle x,h\rangle.
\end{equation}

\item[(ii)] If for some $0\neq h^*\in\real^d$ we have $M_{\pr}(h^*)<\infty$, then for all $x$ such that $\langle x,h^*\rangle > 0$, we have
\begin{equation}\label{eqn:light-tail-1}
    \liminf_{t\to\infty} \frac1t \log \left( \frac1{HD(tx,\pr)}\right) > \langle x,h^*\rangle.
\end{equation}

\item[(iii)] Moreover, under the assumption of $(ii)$, let 
$\{f_i\}_{i=1}^d$ be any orthonormal basis of $\real^d$, then there exists $j\in \{1,\ldots,d\}$ such that
\begin{equation}\label{eqn:one-dir-light}
    \liminf_{t\to\infty} \frac1t \max\left(\log\left( \frac1{HD(tf_j,\pr)}\right),\,\, \log\left( \frac1{HD(-tf_j,\pr)}\right)\right) > |\langle h^*,f_j\rangle|. 
\end{equation}
\end{itemize}
\end{theorem}

The result in Theorem~\ref{thm:hd-decay-light-tail}$(iii)$ is interesting as it implies that, if the measure $\pr$ has a light tail along any direction, then the half-space depth computed along any $d$ directions forming an orthonormal basis, will capture the light tail behaviour with certainty. 

The following proposition identifies various tail behaviours, classifying them in order of increasing heaviness, from super-exponential through exponential and sub-exponential to polynomial. This will be a useful tool for interpreting the HD-plots in terms of tail behaviour when comparing them to the exponential decay ($e^{-t}$), as discussed in Section~\ref{sec:HDmethod}.
\begin{proposition}
\label{prop:tool-variousTails}~
    \begin{enumerate}
    \item[(i)] If $HD(tf_j,\pr) \sim t^{\gamma} e^{-\alpha t^{\beta}}$ for some $\alpha>0,\,\beta> 0,\,\gamma\ge 0$. Then, 
    \begin{eqnarray*}\lim_{t\to\infty} \frac1t \log\left(\frac1{HD(tf_j,\pr)}\right) &=& \lim_{t\to\infty} \left(\alpha \, t^{\beta-1} - \frac{\gamma\log t}{t} \right)\\
    &=&    
    \left\{  
    \begin{array}{cl} 
    +\infty & \text{if } \beta > 1,\\
    \alpha, & \text{if } \beta=1,\\
    0 &  \text{if } \beta<1.
    \end{array}\right. \end{eqnarray*} 
    \item[(ii)] If $HD(tf_j,\pr) \sim (\log t)^{\delta}\,t^{-\theta}$ for some $\delta,\theta>0$. Then, $$\lim_{t\to\infty} \frac1t \log\left( \frac1{HD(tf_j,\pr)}\right) = \lim_{t\to\infty} \frac1t\left(\theta\log t - \delta \log\log t \right)= 0.
    $$
    \item[(iii)] Given that $0$ can be the limit for both, subexponential and polynomial, types of behaviours of HD, we compare $\log\left(1/HD(tf_j,\pr)\right)$ with $\log t$, rather than with $t$, to discriminate between the two behaviours. 
    \end{enumerate}
\end{proposition}

\begin{remark}\label{rk:lower-upper-HD}~
\begin{enumerate}
    \item
    The last result in Theorem~\ref{thm:hd-decay-light-tail} (iii) is interesting as it implies that, if the measure $\pr$ has a light tail along any direction, then the half-space depth computed along any $d$ directions forming an orthonormal basis, will capture the light tail behaviour with certainty. We shall use this precise feature in developing an exploratory tool for multivariate data analysis to identify the tail behaviour of the underlying distribution in Section~\ref{sec:HDmethod}. 
    \item 
    Under the assumptions of Theorem~\ref{thm:hd-decay-tail}, combining the lower bound given in Theorem~\ref{thm:hd-decay-tail} with the upper bound given in Proposition~\ref{thm:hd-decay-survivalMargins},  
we have, for sufficiently large $t$,
\begin{equation*}
c\,g(t\|x\|) \le HD(t x,\pr) \le 
\min\!\Big(\!\min_{1\le i\le d} F_{e_i}(t\,\langle x,e_i\rangle), \min_{1\le i\le d}\left(1-F_{e_i}(t \, \langle x,e_i\rangle)\right)\!\Big),
\end{equation*}
for any $x\ne 0$ and any canonical basis $\{e_i\}_i$. 
\end{enumerate}
\end{remark}

Based on the results obtained, we conclude that the asymptotic behaviour of halfspace depth accurately reflects the asymptotic behaviour of the underlying probability measure. 
Moving from population to sample versions, we will explore whether depth functions effectively mirror the tail behaviour of the underlying distribution; this question will be central to our investigation in Section~\ref{sec:empirical-HD}. 

\section{Empirical multivariate halfspace depth}
\label{sec:empirical-HD}

As noted earlier in the introduction, our motivation for studying halfspace depths is that they are known to capture certain behaviour of the underlying distribution, although they do not uniquely characterise it, as shown by \cite{Nagy2021a}. However, in this last quoted paper, the author listed eight situations for which the halfspace depths uniquely identify the underlying distribution, empirical measure (which is our focus) being one among them. 
Note also that \cite{Struyf1999} proved that halfspace depth uniquely identifies measures with finite support, such as empirical measures.

When analysing the decay of halfspace depth for empirical measures, it is natural to compare their decay rate with those of the underlying measures from which the samples are drawn. This comparison is evaluated in Theorem~\ref{thm:empHDdecay-vs-popHD}.

In the case of heavy tails, a related question has been addressed in \cite{Einmahl2015}, which deals with the challenge of empirical half-space depth vanishing outside the convex hull of the data. In that work, the authors use extreme value statistics to complement the empirical estimator of half-space depth with an additional estimator based on the limiting extreme distribution of half-space depth in the tail, thereby providing a way to estimate depth outside the data hull.

In our setting, while remaining within the data hull, we present a general result on the empirical half-space depth that imposes no assumption on the continuity of the underlying distribution, and covers both light and heavy tails. Moreover, our result is robust, as the convergence holds almost surely.
\begin{theorem}
    \label{thm:empHDdecay-vs-popHD}
    Let $\pr$ be a probability measure defined on $(\real^d,\calB(\real^d))$, and let $g_c$ be the capacity function corresponding to $\pr$, which definition is recalled in Appendix~\ref{subsec:app2}.
    Consider an i.i.d. sample $\{X_k\}_{k\ge 1}$ drawn from $\pr$, and $HD(\cdot,\pr_n)$ the halfspace depth for the empirical measure $\pr_n = \frac1n \sum_{i=1}^n \delta_{X_i}$. 
    
    Let $\{\gamma_n\}_{n\ge 1}$ be a deterministic sequence satisfying the following conditions:
    $$
    \centering
            (C1a) \qquad n^{-1}\log (g_c(\gamma_n)) = {\textrm o}(\gamma_n), 
    $$ 
    and 
    $$
    \centering
            (C1b) \qquad n^{-1}\log \log n = {\textrm o}(\gamma_n). 
    $$
Additionally, \\[-2ex]
    \begin{center}
    $(C2) \qquad$ let $\{t_n\}_{n\ge 1}$ \text{with} $\displaystyle t_n \underset{n\to\infty}{\to}\infty$, such that for some $\veps > 0$, 
    \end{center}
    $$
    \quad HD(t_n x,\pr)>\gamma_n \quad \text{for large enough}\;n \text{ and all } x \in \real^d \text{ with } \|x\| \leq \veps.
    $$ 
    Then, we have
    \begin{equation}\label{eqn:ratio-emp-pop-HD}
        \sup_{\|x\|\leq\veps}\left|\dfrac{HD(t_n x, \pr_n)}{HD(t_n x, \pr)}-1\right| \underset{n\to\infty}{\to} 0\quad \text{almost surely}.
    \end{equation}
    If $\{\gamma_n\}$ satisfies $(C1a)$ only, then, under $(C2)$, \eqref{eqn:ratio-emp-pop-HD} holds in probability.
    \end{theorem}
    \begin{remark}~ 
    \vspace{-1ex}
        \begin{itemize}
        \item[-] Note that $(C1b)$, sufficient condition given in \cite[Theorem~5.1]{Alexander1987} on which we build our result, implies that $\displaystyle n\gamma_n {\to} \infty $ as $n\to\infty$. 
        \\[-4ex]
        \item[-] The capacity function $g_c$ is introduced to define the condition $(C1a)$ in Alexander's result on the convergence of the empirical measure to the underlying parent distribution, which holds regardless of the tail behaviour. We will observe that, although $g_c$ varies depending on the tail behaviour, its influence on the choice of $\gamma_n$ - which satisfies $(C1b)$ -  in order to satisfy $(C1a)$ is negligible. The impact of the tail behaviour will instead be reflected in the empirical halfspace depth through the rate of decay of the halfspace depth of the parent measure $\pr$.
        \\[-4ex]
        \item[-] Condition $(C2)$ implies that, in order to apply Theorem~\ref{thm:empHDdecay-vs-popHD} in any setting, we must have a reasonable way of estimating $HD(t_n x,\pr)$; this is why we have provided estimates for halfspace depths in Section~\ref{sec:multi-quantiles-HD}. Furthermore, it is important to note that the depths are computed at $t_n x$, where $\{t_n\}$ grows with the sample size. In Examples \ref{exple:mrv} and \ref{exple:exp-normal}, we consider different distributions and provide estimates of $\{t_n\}$ that satisfy Condition $(C2)$. Finally, considering $x$ within a bounded set is essential; otherwise, no sequence $\{t_n\}$ could satisfy Condition $(C2)$.
        \\[-4ex]
        \item[-] For convergence in probability, which holds under $(C1a)$ and $(C2)$, a similar result is presented in \cite{He2017}, where the supremum is taken over a central depth region. In contrast, for simplicity, we take the supremum in \eqref{eqn:ratio-emp-pop-HD} over a fixed  $\veps$-ball.
        \\[-4ex]
        \item[-] Finally, as can be seen in the proof of Theorem~\ref{thm:empHDdecay-vs-popHD}, the condition of the halfspace depth is transferred to an appropriate decay condition on the tail probabilities of the measure $\pr$.
        \end{itemize}
    \end{remark}
    We can also propose in a straightforward way the following bound for 
the expression of interest given in \eqref{eqn:ratio-emp-pop-HD}, which will be useful when developing the proof of Theorem~\ref{thm:empHDdecay-vs-popHD}: 
\begin{lemma}\label{prop:hdepth-sup-emp-ratio-gl}
For any $x \in \real^d$, we have the following inequality:
       $$ 
        \left|\dfrac{HD(t_n x, \pr_n)}{HD(t_n x, \pr)}-1\right| \le \sup_{H \in \calH_{t_nx}} \left|\frac{\pr_n(H)}{\pr(H)} - 1 \right|.
       $$  
\end{lemma}

The bound thus obtained has been analysed by many researchers in one or the other form, 
the Glivenko--Cantelli theorem being one of the earliest in this direction. Later, the rate of convergence was obtained by several authors, e.g. \cite{Alexander1987,Gine2006,Shorack2009,Wellner1992,Burr2017} in different scenarios with specific assumptions. It is noteworthy that most of the results in this direction use the specific structure of $\calH_x$ and Dvoretzky--Kiefer--Wolfowitz (DKW) inequality (\cite{Dvoretzky1953}). The specific structure we refer to is called the Vapnik--Chervonenkis (VC) class (see \cite{Vapnik1971,Dudley1984,Alexander1987,Talagrand2003}). The idea of VC class has its roots in statistical learning wherein one is interested in identifying the class of functions to characterise convergence of probability measures. Specifically, a class $\calS$ of sets shatters a finite set $F$ if, given $G\subset F$, $\exists S\in\calS$ for which $G = F \cap S$. A class $\calS$ of sets is called a VC class if for some integer $n$, $\calS$ shatters a set of cardinality $n$ and does not shatter any set of cardinality bigger than $n$ (see Appendix~\ref{subsec:app2} for more details).
In our analysis, we shall use the approach of \cite{Alexander1987} on VC class, without resorting to the DKW inequality. Specifically, we shall invoke
Theorem~5.1 of the \cite{Alexander1987}, which we recall in Appendix~\ref{subsec:app2}, powerful and crucial result for the proof of Theorem~\ref{thm:empHDdecay-vs-popHD}. 

Notice that, without any assumption on $\pr$, Theorem~\ref{thm:empHDdecay-vs-popHD} presents a general result about the rate of decay of the empirical halfspace depth and allows us to compare the halfspace depth of the parent measure $\pr$ and of the empirical measure $\pr_n$. Indeed, \eqref{eqn:ratio-emp-pop-HD} implies that
\begin{equation}\label{eqn:emp-pop-hdepth}
    \pr\Big(c_1\, HD(t_nx, \pr)\,\le \, HD(t_n x, \pr_n)  \,\le  \,c_2 HD(t_nx, \pr)\Big) \to 1\,\,\,\,\text{ as } n\to\infty
\end{equation}
for all $c_1 \le 1\ \le c_2$. We will have more to say about it when we specialise to specific cases.

However, as seen in Section~\ref{sec:multi-quantiles-HD}, the rate of decay of halfspace depth is closely related to the tail behaviour of $\pr$, which, in view of \eqref{eqn:emp-pop-hdepth}, implies that the rate of decay of the empirical halfspace depth can be estimated as a function of the tail behaviour of the parent measure $\pr$. The following three theorems establish this connection.

We start with a lower bound deduced from Theorem~\ref{thm:hd-decay-tail}.
\begin{theorem}
\label{thm:empirical-hd-decay-tail}
    If there exists a positive function $g$ satisfying 
    the condition \eqref{eqn:unif-direction-density-decay} given in Theorem~\ref{thm:hd-decay-tail}, and if $(\gamma_n)$ satisfy Conditions $(C1a)$, $(C1b)$ and $(C2)$ in Theorem~\ref{thm:empHDdecay-vs-popHD}, then there exists $c>0$ such that
      \begin{center}$\displaystyle
            \pr\left(\liminf_n \frac{HD(t_nx,\pr_n)}{g(t_n\|x\|)} \ge c\right) =1 
    $\end{center}
\end{theorem}
This implies that, for sufficiently large $n$, the empirical halfspace depth can be lower bounded by the decay function specified in \eqref{eqn:unif-direction-density-decay} with high probability.

Under the similar MRV framework as in \cite[Proposition~2]{He2017}, recalled in Proposition~\ref{prop:depth_decay_pop_nodensity}, we obtain in the heavy tail case the following rate of convergence for the empirical halfspace depth:

{\begin{theorem}\label{thm:emp-hdepth-mrv-decay_nodensity}
    Let $\pr$ be a probability measure defined on $(\real^d,\calB(\real^d))$ that satisfies the following conditions: 
    \begin{enumerate}
        \item There exists some measure $\nu$ such that
        \begin{equation*}
        \lim_{t \to \infty}\dfrac{\pr(tA)}{1-\pr(tB^d)}= \nu(A)<\infty
        \end{equation*}
    for every Borel set $A \subset \real^d$ that is bounded away from origin and satisfies $\nu(\partial A)=0$, where $tA=\{t x: x \in A\}$. Additionally, $\nu(A)>0$ if $A\supset H$ for some halfspace $H$.
        \item $\pr( \calC_{\beta} )=0 $ where $\calC_{\beta} = \{ x : HD(x,\pr)= \beta\}$ is the depth contour at level $\beta \in (0,1)$.
        \item $\displaystyle \lim_{t\to\infty}\dfrac{1-\pr(tB^d)}{t^{-\alpha}}=c \in (0, \infty)$ for some $\alpha>0$.
    \end{enumerate}    
    Let $\{X_n\}_{n\ge 1}$ be an i.i.d sample drawn from $\pr$. Then, we have,
    \begin{equation}
        \lim_{n \to \infty} \sup_{{\|x\|=1} }\left|\frac{HD(t_n x,\pr_n)}{1-\pr(t_n B^d)} - HD(x, \nu) \right| = 0, \,\,\,\,\,\,\,\text{a.s.}
    \end{equation}
whenever $(t_n)\nearrow \infty$ such that $1-\pr(t_n B^d)>\, \gamma_n$ (for any large $n$), and $(\gamma_n)$ satisfies Condition  $(C1b)$ given in Theorem \ref{thm:empHDdecay-vs-popHD}. 
\end{theorem}
}

\begin{remark}~
\vspace{-1ex}
    \begin{itemize}
    \item[-]  Let $g_c$ be the capacity function corresponding to $\pr$. Since $g_c(\cdot)=O(1)$ in the MRV setting (see Appendix~\ref{subsec:app2}), Condition $(C1b)$ implies Condition $(C1a)$. Consequently,
    $(C1a)$ is equivalent to the condition: $n\gamma_n \to \infty$ as $n\to \infty$.
    \\[-4ex]
    \item[-] If assuming Condition $(C1a)$ only, the result holds in probability only, as given in \cite{Einmahl2015}, Equation~(31) and in \cite[Proposition 2]{He2017}.
    \\[-4ex]
    \item[-] Further, we note that the focus in \cite{Einmahl2015} is on extrapolating the empirical depth beyond the data cloud in the MRV setting (subclass of heavy tail distributions), which allows them to relax Condition $(C1a)$ into $n\gamma_n \to 0$ (instead of $\infty$) as $n\to \infty$. As the authors use extreme value statistics, additional conditions on $\gamma_n$ are required, which are common in extreme value theory. The focus of our study is different, which is why we do not consider such extrapolation. Instead, we adopt a more general framework that emcompasses both light and heavy tails and allows us to obtain almost sure convergence.
    \\[-4ex]
    \item[-] Theorem~\ref{thm:emp-hdepth-mrv-decay_nodensity} can be rewritten with explicit conditions when the probability measure has a regularly varying density; we present the result in Appendix~\ref{app:He2017} for completeness, to recall the explicit conditions given on the density (see Corollary~\ref{cor:theo-emp-hdepth-mrv}).
    \end{itemize}
\end{remark}

\begin{example}\label{exple:mrv}~
    For $\alpha>0$, consider a MRV distribution with index $-\alpha$, then $1-\pr(t_n B^d)\approx t_n^{-\alpha}$. Therefore, Theorem~\ref{thm:emp-hdepth-mrv-decay_nodensity} holds if $t_n\nearrow \infty$ and $t_n \leq \gamma_n^{-1/\alpha}$. Now, by choosing $\gamma_n:=n^{-\beta}$ with $0<\beta<1$, the condition $1-\pr(t_n B^d)>\, \gamma_n$ gives a speed of $t_n \leq n^{\beta/\alpha}$.
We could slightly improve the rate by considering $\gamma_n:=\log_p(n)/n$ with $p> 1$. 
In fact, if targeting the convergence in probability (instead of a.s.), we could choose $\gamma_n:=k_n/n$ with  $k_n \to \infty$ and $\gamma_n\to 0$ as $n\to\infty$.
\end{example}

Let us turn to the light tail case for which we can provide, under distinct conditions, a lower and an upper bound for the asymptotics of the halfspace depth. Through examples, we observe that the general bounds can be tight as in the exponential case (see Example~\ref{exple:exp-normal},(a)), with the lower bound of the order of the upper one, showing that this rate is optimal. In the case of spherically symmetric distributions, the gap between the bounds can be reduced via a direct computation as given in Example~\ref{exple:exp-normal},(b) below.

\begin{theorem}\label{thm:emp-hdepth-light-decay}
Let $\{X_n\}_{n\ge 1}$ be an i.i.d. sample drawn from a probability measure $\pr$ defined on $(\real^d,\calB(\real^d))$.
Assume that the moment generating function 
$M_{\pr}(h^*)<\infty$  for some $0\neq h^*\in\real^d$.
Choosing $(\gamma_n)$ satisfying Conditions $(C1a)$, $(C1b)$ and $(C2)$ in Theorem~\ref{thm:empHDdecay-vs-popHD}, we have, for $x$ such that $\langle x,h^*\rangle > 0$,
    \begin{equation}\label{eqn:emp-light-tail}
     \pr\left[\liminf_n \frac{1}{t_n} \log \left( \frac1{HD(t_nx,\pr_n)}\right)  > \langle x,h^*\rangle \right] =1 
\end{equation}

\end{theorem}

\begin{remark}~
    \begin{enumerate}
        
        \item Assuming the conditions of Theorems~\ref{thm:empirical-hd-decay-tail}\,\&\,\ref{thm:emp-hdepth-light-decay} and choosing $t_n$ such that $g(t_n)> \gamma_n$ ($g$ being introduced in Theorem~\ref{thm:empirical-hd-decay-tail}), we can write 
    \begin{equation}\label{eqn:hd-lower-upper-light}
    c\,\gamma_n \le HD(t_n x,\pr_n) \le C\,e^{-t_n\langle x,h^*\rangle},
    \end{equation}
    with high probability for large enough $n$, for some finite, positive $c$ and $C$, that do not depend on $n$. 
    \item Using Proposition~\ref{thm:hd-decay-survivalMargins} for the empirical measure $\pr_n$ and the relation \eqref{eqn:emp-pop-hdepth} between halfspace depth of $\pr$ and of $\pr_n$, the upper bound given in \eqref{eqn:hd-lower-upper-light} can further be refined. Specifically, writing $F^{(n)}_i$ for the $i$-th empirical marginal distribution of $\pr_n$, we have
    \begin{eqnarray}\label{eq:upperBd-improved}
    HD(t_n x,\pr_n) &\le & \bar{C}\,\min \Big(e^{-t_n \langle x,h^*\rangle}, \\
    && \qquad \min_{1\le i\le d} (1-F^{(n)}_i(t_n \, x_i)), \min_{1\le i\le d} F^{(n)}_i(t_n \, x_i)\Big), \nonumber
    \end{eqnarray} 
    with high probability for large enough $n$, for some finite positive $\bar{C}$.
    \end{enumerate}
\end{remark}

\begin{example}\label{exple:exp-normal}~     
    \begin{itemize}
        \item[(a)] {\bf Exponential case.} 
        Let $Y$ be a random vector with probability density function $h$ defined by
    $\displaystyle h(y)=ke^{-\|y\|}$ for $y \in \real^d$ and $k>0$, the normalising constant.
    Then, we have     
\begin{equation}\label{eq:ExponExple}
    \lim_{R \to \infty} \inf_{\theta\in S^{d-1}}\dfrac{\pr\big( \langle\theta, Y\rangle \geq R\big)}{e^{-R}} >0,
\end{equation}
from which we deduce that $c_1\, e^{-t_n} \le HD(t_nx,\pr_n) \le c_2 \,e^{-t_n}$ with high probability for large enough $n$ and for some finite $c_1, c_2>0$.\\[1ex]
\textnormal{In order to prove \eqref{eq:ExponExple}, observe that, due to the spherical symmetry of the distribution of $Y$, for any $\theta_1, \theta_2\in S^{d-1}$, 
$\displaystyle
\pr\big( \langle\theta_1, Y\rangle \geq R\big)=\pr\big( \langle\theta_2, Y\rangle \geq R\big)$.
Therefore, we can consider any direction and, choosing $\theta=(1,0,...,0)$ and using Minkowski inequality, we can write,
\begin{align*}
        &\inf_{\theta\in S^{d-1}}\pr\big( \langle\theta, Y\rangle \geq R\big)
        =\pr\big( Y_1 \geq R\big)
        =\int_{-\infty}^{\infty}\cdots \int_{-\infty}^{\infty}\int_R^{\infty} ke^{-\|y\|}dy\\
        & \qquad\geq \int_{-\infty}^{\infty}\cdot \cdot \cdot \int_R^{\infty} ke^{-(|y_1|+...+|y_d|)}dy =K \int_R^{\infty} e^{- y_1} dy_1
        =Ke^{-R},
\end{align*}
for some constant $K>0$.\\
We deduce that \;$\displaystyle \frac{\inf_{\theta\in S^{d-1}}\pr\big( \langle\theta, Y\rangle \geq R\big)}{e^{-R}}\ne 0$ for large enough $R$, hence \eqref{eq:ExponExple}.} \hfill \qed
        \item[(b)] {\bf Gaussian case.} Consider a multivariate standard normal distribution.
        We choose $g(R)=e^{-R^2 / 2}$ (based on the Mill's ratio). 
        Now, set $\gamma_n=n^{-\beta}$ with $0<\beta <1$. 
        Such $\gamma_n$ satisfies Conditions $(C1a)$ and $(C1b)$ by way of Example 3.6 of \cite{Alexander1987} (see also Appendix~\ref{subsec:app2}). 
        The condition $g(t_n)> \gamma_n$ gives $t_n \leq \sqrt{2 \beta \log n}$, and, for such $t_n$, 
        $$
        HD(t_n x,\pr_n)> c\,n^{-\beta}
        $$
        with high probability for large enough $n$.
        The upper bound as given in \eqref{eqn:hd-lower-upper-light}, of order $e^{-t_n}$, is then
        quite large compared with the lower bound. Nevertheless, it can be improved, simply by considering the Gaussian marginal distributions as in \eqref{eq:upperBd-improved}. In such a case, the upper bound for $\displaystyle t_n = \sqrt{2 \beta \log n}$ becomes
        $HD(t_n x,\pr_n)\le C\,n^{-\beta}$ with high probability for large enough $n$.
    \end{itemize}
\end{example}

\section{An exploratory HD-tool to discriminate between light and heavy tails}
\label{sec:HDmethod}

\subsection{A methodology to assess the tail behaviour of a given multivariate dataset}
\label{ss:HDtool}

Based on the theoretical results obtained above, we shall, in this section, derive a methodology to unravel the tail behaviour of a given multivariate dataset. We present the methodology as a pseudo-algorithm given below (Algorithm \ref{algo}).

\begin{algorithm}
\caption{Tail Behavior Detection via Halfspace Depth (Concise)}
\label{algo}
\begin{algorithmic}[1]
\FOR{$k = 1$ to $d$}
    \STATE Compute and plot $\displaystyle y_n = \frac{1}{t_n} \log\left( \frac{1}{HD(t_n e_k, \pr_n)} \right)$ vs $t_n$ for $n = 1, 2, \dots$
    \STATE\textit{Visually analyze the asymptotic behavior of $y_n$ as $t_n \to \infty$}
    \IF{$y_n \to \infty$ or $y_n \to c > 0$}
        \STATE Flag direction $k$ as \textbf{light-tailed}
    \ELSE
        \STATE Compute and plot $\displaystyle w_n = \frac{1}{\log(t_n)} \log\left( \frac{1}{HD(t_n e_k, \pr_n)} \right)$ vs $t_n$
        \IF{$w_n \to \infty$}
            \STATE Flag direction $k$ as \textbf{light-tailed}
        \ELSIF{$w_n \to \theta \in (0, \infty)$}
            \STATE Flag direction $k$ as \textbf{heavy-tailed}
        \ENDIF
    \ENDIF
\ENDFOR
\IF{All directions are flagged as heavy-tailed}
    \STATE \textbf{Conclude}: The distribution is heavy-tailed
\ELSE
    \STATE \textbf{Conclude}: Light-tailed along some direction
\ENDIF
\end{algorithmic}
\end{algorithm}

If the lightest tail discovered through the algorithm corresponds to a light tail (exponential type), then, by Theorem~\ref{thm:emp-hdepth-light-decay}, we can deduce that there is strong evidence that $\pr$ has a light tailed marginal distribution along some direction. In such a case, $HD(x,\pr_n)$ provides a good estimate for $HD(x,\pr)$ (with no need of extrapolation beyond the convex hull of the sample considered). On the other hand, if the lightest $HD(t_n x,\pr_n)$ corresponds to a heavy tail, then, by Theorem~\ref{thm:emp-hdepth-mrv-decay_nodensity}, we can deduce that $\pr$ has a heavy tail. In this case, the refined estimator defined in  \cite{Einmahl2015} by extrapolating beyond the hull can be used.

\subsection{Application on simulated data and discussion}
\label{ss:applic}

Consider a 3-dimensional random vector ${\bf X}$ and $A \in SO(3)$ such that $A\neq I_{3\times 3}$, such that ${\bf Y} = A {\bf X}$ has independent marginals with $Y_1\sim N(0,1)$, $Y_2\sim \text{Laplace}$, and $Y_3\sim t_3$, the Student's $t$-distribution with parameter $3$.

We generate $N$ i.i.d. copies (samples) ${\bf X_1},\ldots,{\bf X}_N$ of ${\bf X}$. Writing $(e_1,e_2,e_3)$ for the canonical basis of $\real^3$, we set $X_{i,j} = \langle e_i,{\bf X_j}\rangle$ as the projection of ${\bf X_j}$ onto the $e_i$ direction.

For a fixed positive integer $M$, we set $n=\lfloor kN/M\rfloor$, for $k=1,\ldots,M$. Writing $\pr_n$ for the empirical measure corresponding to the sample $\{{\bf X}_j\}_{j=1}^n$, we compute $HD(t_ne_i,\pr_n)$, where $t_n = n/{1000}$ for different values of $n$ and $i$.

Notice that, if $A_{3,i}\neq 0$, for all $i,j=1,2,3$, then the distribution of the $i$-th marginal $X_i$ will exhibit a heavy tail behaviour due to the contribution of $Y_3$, for all $i=1,2,3$. Note that the effect of $Y_3$ on $X_i$ will also depend on the sign of $A_{3,i}$. Specifically, if $A_{3,i}>0$, then the tails of $Y_3$ and $X_i$ align. 
Conversely, if $A_{3,i}<0$, then the right tail of $Y_3$ will impact the left tail of $X_i$, and left tail of $Y_3$ will impact the right tail of $X_i$. Nevertheless, in our case, $Y_i$ have symmetric distribution, hence, the sign of $A_{ij}$ will not impact the tail behaviour of $X_i$s.
\\Now, considering another basis $\{f_i\}_{i=1}^3$ of $\real^3$ defined by $\displaystyle f_i= A^T e_i$, we can write, using \eqref{eq:gl-upperBd-HD} and \eqref{eq:upperBd-improved},
$$
HD(t_n e_i,\pr_n) \le \min \left( \min_{1\le j\le 3}\!\!\left( 1 - F_j^{(n)}\!(t_n \langle e_i,f_j \rangle)\right), \min_{1\le j\le 3} F_j^{(n)}\!(t_n \langle e_i,f_j \rangle)\right),
$$
where $F_j^{(n)}$ denotes the marginal cumulative distribution function of $\pr_n$ along the $f_j$-direction, which is also the empirical marginal of $Y_j$.

\begin{figure}[h]
    \centering
    \begin{minipage}{0.44\textwidth}
       \includegraphics[width=0.9\linewidth]{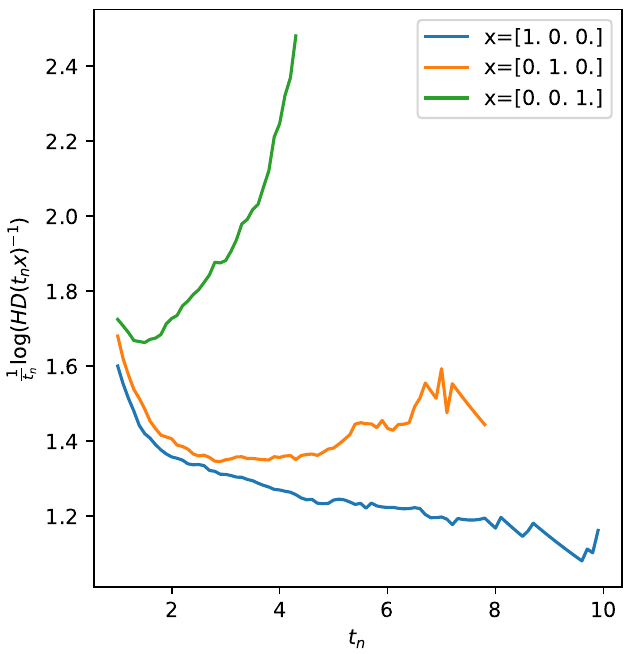}
    \end{minipage}   
    \hfill
    \begin{minipage}{0.44\textwidth}
       \includegraphics[width=0.9\linewidth]{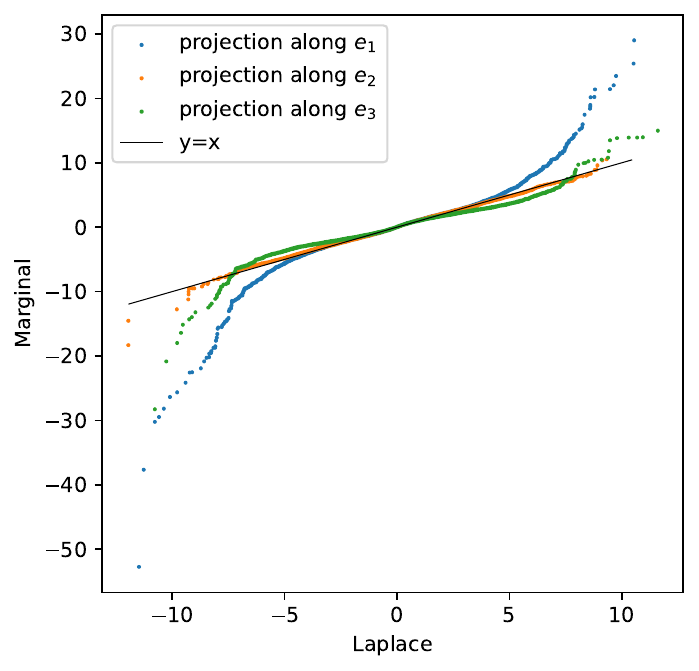}
    \end{minipage}   
\parbox{330pt}{\caption{\sf \small 
Left plot: $\displaystyle y_n := \log\left(\left( {HD(t_n x, \pr_n)} \right)^{-1}\right)/{t_n} $ against $t_n=n/1000$ (with $n=k10^3$, for $k=10,\ldots,100$) in $e_1, e_2$ and $e_3$ directions.
Right plot: QQ-plot of the marginals $X_i$, $i=1,2,3$ (where ${\bf X} = A^T{\bf Y}$), w.r.t. the standard Laplace quantiles.  
\label{fig:simulData-HDrate}}}
\end{figure}

For our example, illustrated in Figure~\ref{fig:simulData-HDrate}, we take $N=10^5$, $M=100$ and 
$$
A=\left[\begin{array}{ccc} 0.3536 & -0.4189 &  0.8364 \\ 0.3536 &  0.8876 &  0.2952 \\
-0.8660 & 0.1913 & 0.4619 \end{array}\right].
$$ 
Notice that we need not consider the $-e_i$ directions, as all the marginals are symmetric. 
The QQ--plot clearly shows that $X_1, X_2$ and $X_3$ exhibit heavy-tailed behaviour. In contrast, the halfspace depth plots in the directions $e_1$, $e_2$ and $e_3$ (left plot) exhibit mixed behaviour. 
For the $e_1$ direction, no definitive conclusion can be drawn: Although the trend appears to be decreasing, we cannot extend $t_n$ further, since all available observations have been used (and the extrapolation method suggested by \cite{He2017} does not help in this case as it relies on the MRV assumption). Therefore, we cannot determine visually its asymptotic behaviour. 
On the other hand, the $e_2$ and $e_3$ directions appear to exhibit exponential and Gaussian--type behaviour, respectively. Note that since, $\langle f_i,e_j\rangle = A_{ij}\neq 0$, it is expected that the Gaussian tail would asymptotically dominate the tail behaviour of $HD(t_n e_i,\pr_n)$. However, we observe an apparent non--Gaussian tail in our experiment along the $e_1$ and $e_2$ directions. At first glance, we attribute this deviation primarily to the relatively small range of $t_n$ in the experiment, which is insufficient to capture the large sample asymptotics. In such cases of relatively smaller values of $t_n$, the coefficient of $t_n$ in $F_j(t_n\langle e_i,f_j\rangle)$ plays a significant part in determining the bound for $HD(t_ne_i,\pr_n)$. Specifically, notice that in case of $HD(t_ne_1,\pr_n)$, the coefficient $\langle f_1,e_3\rangle$ is significantly larger than $\langle f_1,e_1\rangle$ and $\langle f_1,e_2\rangle$. Thus, at the scale at which we can vary $t_n$ in this experiment, the decay of $F_1$ (polynomial) appears faster than the exponential and Gaussian. Similarly, for $HD(t_ne_2,\pr_n)$, the coefficient $\langle f_1,e_2\rangle$ is significantly larger than $\langle f_1,e_1\rangle$ and $\langle f_1,e_3\rangle$, which in turn leads to the exponential function dominating the polynomial and Gaussian. In summary, we conclude that although not all $e_i$ directions could exhibit the expected Gaussian behaviour, at least one direction was able to capture the underlying Gaussian tail.

To conclude the discussion, the following are key uses for the HD tool: 
\\ (i) It serves as an excellent tool for comparing distributions. 
\\ (ii) In practice, the HD tool works well as a complementary method to QQ-plots. While QQ-plots are often insightful, they may not always provide a clear picture. In such cases, the HD tool can help to further clarify the analysis.
\\ (iii) The HD tool can also help assess the heaviness or lightness of the tail, although this is more challenging. Nevertheless, it cannot substitute for statistical tests specifically designed to assess tail behaviour. Note that this was already the case in \cite{He2017}, where HD extrapolation for extremes is discussed in the MRV setting. This means one should first verify whether the data exhibit characteristics of MRV before proceeding with HD-based tail analysis.

\subsection{Application on real data}
\label{ss:appli-RealData}

We consider the data corresponding to {\it outgoing longwave radiation (OLR)} measured at various locations in India. The dataset and its detailed description can be accessed from \url{https://tropmet.res.in/static_pages.php?page_id=144}. At large scale, OLR can be interpreted as a proxy of cloud cover: a relatively high value of OLR indicates clear skies, while small OLR values correspond to cloudy skies. Therefore, in the Monsoon zone in India, the OLR is more likely to take relatively small values during the summer Monsoon months of June, July, August and September, compared with the Winter months of January, February, March and April. Therefore, we expect the OLR distribution to exhibit a heavier left tail during the Summer Monsoon season. In contrast, the OLR values are mostly concentrated around reasonably high values during the Winter season when the skies are mostly clear, thus, indicating a lighter tail.

We assess the OLR distribution at two locations (Jaisalmer and Prayagraj), in two different seasons (the Summer Monsoon period of June, July, August and September and Winter months of January, February, March and April) for the years $2004$ to $2017$. Jaisalmer, a city in North Western part of India is known to experience arid weather throughout the year with annual precipitation of $17$cm and is typically not considered in the Monsoon zone, whereas Prayagraj is located in the Monsoon zone and receives $104$cm of rain annually, mostly during the Monsoon months of June, July, August and September. These two locations are chosen as they present a reasonably contrasting image of monsoon in India. We center the data for better interpretability of the plots that follow.

We consider the $4$ dimensional dataset: Jaisalmer winter, Jaisalmer summer, Prayagraj winter and Prayagraj summer, denoted by $e_1$, $e_2$, $e_3$, and $e_4$ directions, respectively. 
For $N = 1708$, $n=kN/100$, for $k=10,\ldots,100$, we set $\pr_n$ as the empirical measure corresponding to the first $n$ observations in the dataset. 

\begin{figure}[h]
    \centering
    \begin{minipage}{0.32\textwidth}
        \includegraphics[width=0.9\linewidth]{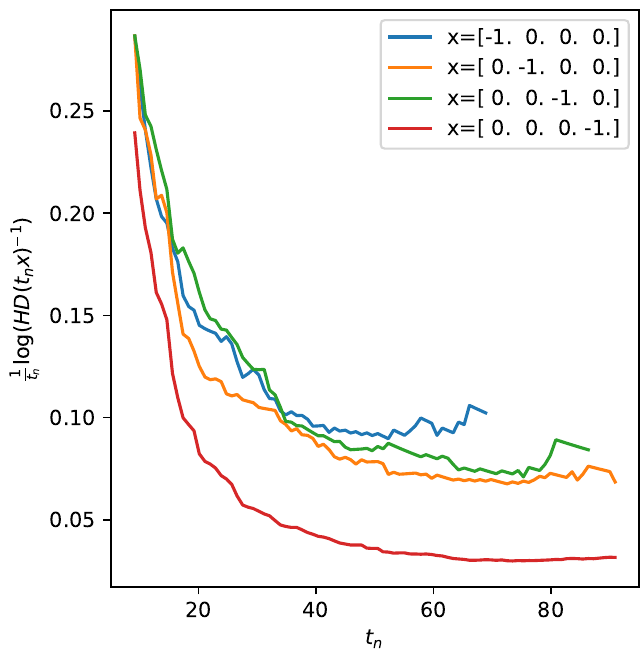}
    \end{minipage}
    \hfill
    \begin{minipage}{0.32\textwidth}
       \includegraphics[width=0.9\linewidth]{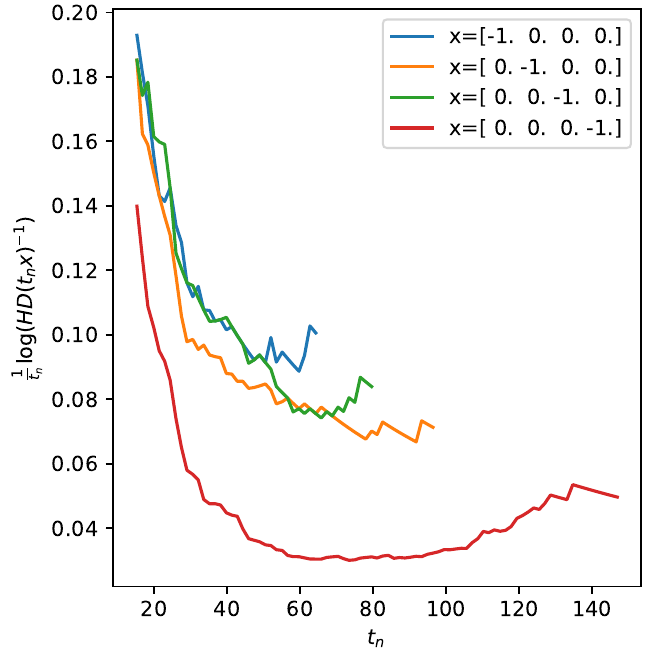}
    \end{minipage}
    \hfill
    \begin{minipage}{0.32\textwidth}
       \includegraphics[width=0.9\linewidth]{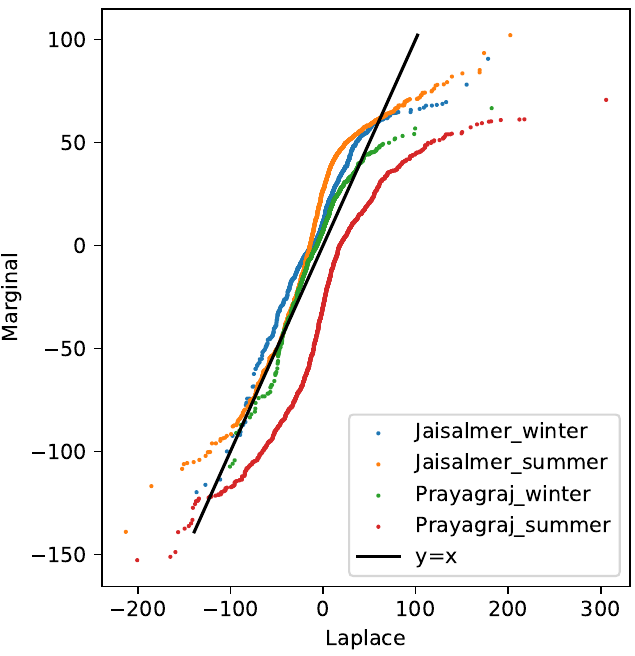}
    \end{minipage}   
\parbox{330pt}{\caption{\sf \small The left and middle plots correspond to $\displaystyle y_n := \log\left(\left( {HD(t_n x, \pr_n)} \right)^{-1}\right)/{t_n} $ plotted for the OLR data in $-e_1, -e_2, -e_3$ and $-e_4$ directions against $t_n$, with $t_n = 6n/100$ for the left plot, and $t_n=n/10$ for the middle plot, for $n=\lfloor 17.08 k\rfloor$, for $k=10,\ldots,100$. The right most figure depicts the QQ plot of all the four marginals w.r.t. univariate standard Laplace distribution.}
\label{fig:realData-HDrate-LeftTail}}
\end{figure}

As climatologists observe homogeneous right-tail behaviour of OLR across seasons but heterogeneous left-tail behaviour, we focus on the left tail to better illustrate our methodology. We, therefore, plot $HD(t_n x,\pr_n)$ for $x=-e_1,-e_2,-e_3,-e_4$ in Figure \ref{fig:realData-HDrate-LeftTail}, but for the sake of completeness, we also display the decay of halfspace depth along $e_1,e_2,e_3,e_4$ directions in Appendix ~\ref{app:OLR}.

Figure \ref{fig:realData-HDrate-LeftTail} reveals a clear contrast in tail behavior, with Prayagraj in summer showing the heaviest tail among the four. Meaning, that OLR is often significantly smaller than its mean at Prayagraj in summer (monsoon), which is indeed explained by the summer monsoon, whereas significantly smaller values (relative to the mean) are less likely to occur at Jaisalmer (summer and winter) and at Prayagraj during winter. Also, observe that the halfspace depth for Prayagraj summer for $t_n \lessapprox 80$ appears to suggest a heavy tail bevaiour, however as $t_n$ grows beyond $80$, an upward trend in $\log\left(\left( {HD(t_n x, \pr_n)} \right)^{-1}\right)/{t_n}$ alludes to an exponential tail with a very small exponent (as compared to the other three tail behaviours).

\section{Proofs}
\label{sec:proofs}

We divide this section into two parts: Section \ref{ss:population-proofs} has the proofs of all the results from Section 2, and whereas all the proofs of results related to the empirical halfspace depths are presented in Section \ref{ss:sample-proofs}.

\subsection{Proofs of Section 2}\label{ss:population-proofs}

Proofs of Propositions \ref{thm:hd-decay-survivalMargins} and \ref{prop:tool-variousTails} are not provided as they can either be deduced directly from known definition, or involve a small algebraic computation.

\subsubsection{Proof of Theorem~\ref{thm:asymp-depth-ell-symm}}

Let us consider $x,y\in\real^d$ such that $\|\Sigma x\| = \|\Sigma y\|$. Let $A$ be the orthogonal matrix such that $\Sigma x = A\Sigma y$. Setting $A^* = \Sigma^{-1} A\Sigma$, we have $x = A^*y$, and $\text{det}(A^*) = 1$. Recalling  the Definition~\ref{def:HD} of halfspace depth, 
and setting $$\displaystyle p_y = \text{arg}\min_{\|p\|\neq 0}\int_{\{w:\langle w-t y,p \rangle\geq 0\}}\!\!\! f(w) dw$$ we have, 
$$ \frac{HD(tx;\pr)}{HD(ty;\pr)} =  \frac{\inf_{\|p\|\neq 0} \int_{\{w: \langle w-tx,p \rangle\geq 0\}}f(w) dw}{\int_{\{w:\langle  w- ty, p_y \rangle \geq 0\}}f(w) dw} \le \frac{\int_{\{w: \langle w-tx,Bp_y \rangle\geq 0\}}f(w) dw}{\int_{\{w:\langle  w- ty, p_y \rangle \geq 0\}}f(w) dw},$$
for any linear transformation $B$.

Now, using the transformation $w=tz$, we obtain
$$ 
\frac{HD(tx;\pr)}{HD(ty;\pr)} \le  \frac{\int_{\{z: \langle z-x,Bp_y \rangle\geq 0\}}f(tz) dz}{\int_{\{z:\langle  z- y, p_y \rangle \geq 0\}}f(tz) dz}.
$$
Recalling that $x=A^*y$, we introduce the transformation $z=A^*u$ in the numerator to observe that 
$$
\frac{HD(tx;\pr)}{HD(ty;\pr)} \le  \frac{\int_{\{u: \langle A^*(u-y),Bp_y \rangle\geq 0\}}f(tA^*u) du}{\int_{\{z:\langle  z- y, p_y \rangle \geq 0\}}f(tz) dz}.
$$
Since the above inequality is satisfied for any nonsingular $B$, we set $B^T=(A^*)^{-1}$, to obtain
$$ 
\frac{HD(tx;\pr)}{HD(ty;\pr)} \le  \frac{\int_{\{u: \langle u-y,p_y \rangle\geq 0\}}f(tA^*u) du}{\int_{\{z:\langle  z- y, p_y \rangle \geq 0\}}f(tz) dz}.
$$ 
Now, while observing that $\|\Sigma A^*u\| = \|\Sigma u\|$, we invoke the assumption of asymptotic elliptical symmetry of $f$, to conclude that there exists $t_0$ large enough, such that
$$ \frac{HD(tx;\pr)}{HD(ty;\pr)} \le  \frac{\int_{\{u: \langle u-y,p_y \rangle\geq 0\}} \left[\frac{f(tA^*u)}{f(tu)}\right] f(tu) du}{\int_{\{z:\langle  z- y, p_y \rangle \geq 0\}}f(tz) dz} \le (1+\epsilon),\,\,\,\forall t\ge t_0$$

Using similar arguments, we can also conclude that
$$
\frac{HD(tx;\pr)}{HD(ty;\pr)} \ge (1-\eps), \,\,\,\,\,\,\,\,\forall t\ge t_0,
$$
which concludes the result. \hfill \qed

\subsubsection{Proof of Theorem~\ref{thm:hd-decay-tail}}
\label{ss:theo2.12prev}

Recall that (\eqref{eqn:hdepth-project}), 
\begin{eqnarray*}
HD(x,\pr) &=& \inf_{\theta \in S^{d-1}} \pr\left[\langle \theta,X\rangle \ge \langle \theta,x\rangle\right]\\
&=& \inf_{\substack{\|\theta\|=1\\ \langle\theta,x\rangle \ge 0}}\,\,\,\min \Big(\pr\left[\langle \theta,X\rangle \ge \langle \theta,x\rangle\right], \pr\left[\langle \theta,X\rangle \le \langle \theta,x\rangle\right]\Big)
\end{eqnarray*}
Let $c>0$ be such that, for all sufficiently large, positive $R$,
\begin{equation}
\inf_{\theta \in S^{d-1}}\dfrac{\pr(\langle \theta, X\rangle \geq R)}{g(R)} >c>0.
\end{equation}
Note that, $$\inf_{\theta \in S^{d-1}} \pr(\langle \theta, X\rangle \geq R) = \inf_{\theta \in S^{d-1}} \,\,\min \Big(\pr(\langle \theta, X\rangle \geq R), \,\, \pr(\langle \theta, X\rangle \leq -R)\Big)$$

Combining the definition of halfspace depth and the above lower bound, we have
\begin{eqnarray*}
HD(t x, \pr) &=&  \inf_{\substack{\|\theta\|=1\\ \langle\theta,x\rangle \ge 0}}\,\,\,\min \Big(\pr\left[\langle \theta,X\rangle \ge \langle \theta,x\rangle\right], \pr\left[\langle \theta,X\rangle \le \langle \theta,x\rangle\right]\Big) \\
&\geq & \inf_{\theta \in S^{d-1}} \,\,\min \Big(\pr(\langle \theta, X\rangle \geq t \|x\|), \,\, \pr(\langle \theta, X\rangle \leq -t\|x\|)\Big) \\
&\geq &  c\, g(t\|x\|),
\end{eqnarray*}
hence the lower bound. \hfill \qed

\subsubsection{Proof of Theorem~\ref{thm:hd-decay-light-tail}}
\label{ss:theo2.12}

$(i)$ The proof of the upper bound is based on a simple application of Markov inequality. Let $\calN$ be such that $\E(e^{\langle h,Y\rangle}) <\infty$ for any $Y$ with distribution $\pr$. Then,
\begin{eqnarray*}
HD(tx,\pr) &\le & \inf_{\substack{\|h\|=1 \\ \langle h,x\rangle \ge 0}} \,\, \pr(\langle Y,h\rangle \ge t \langle x,h\rangle) \\
&=& \inf_{h \in \calN} \,\, \pr(\langle Y,h\rangle \ge t \langle x,h\rangle) \\
&= & \inf_{h \in \calN} \pr(e^{\langle Y,h\rangle }\ge e^{t \langle x,h\rangle })
\le \inf_{h \in \calN} e^{-t\langle x,h\rangle }\;\E[e^{\langle Y,h\rangle}] \\
&\le & C_{\calN} \inf_{h\in\calN} e^{-t\langle x,h\rangle},
\end{eqnarray*}
where $C_{\calN} = \sup_{h\in\calN} \E[e^{\langle Y,h\rangle}]$. Implying that,
$$\log\left( \frac1{HD(tx,\pr)}\right) \ge -\log C_{\calN} + t\,\sup_{h\in\calN}\langle h,x\rangle$$
from which the result follows.\\

$(ii)$ Observe that,
$$
HD(tx,\pr) 
\le \pr[\langle Y,h^*\rangle \ge t \langle x,h^*\rangle] = \pr[e^{\langle Y,h^*\rangle }\ge e^{t \langle x,h^*\rangle }]
\le e^{-t\langle x,h^*\rangle }\;\E[e^{\langle Y,h^*\rangle}],
$$
which leads to the conclusion.\\

$(iii)$ Beginning with $h^*$, let us construct an orthonormal basis $\{h^*,h_2,\ldots,h_d\}$. Let $\{f_i\}$ be any other orthonormal basis. We begin with the claim that, there exists $j\in\{1,\ldots,d\}$ such that $\langle f_j,h^*\rangle \neq 0$. Let us assume to the contrary that $\langle f_j,h^*\rangle =0$ for all $j=1,\ldots,d$. Since, $\{f_j\}$ is assumed to be an orthonormal basis of $\real^d$, therefore, any element of $\real^d$ must be expressible as a linear combination of $\{f_j\}$, implying that
$h^* = \sum_{j=1}^d \langle h^*,f_j\rangle f_j$. Combining this with our assumption, that $\langle h^*,f_j\rangle=0$ for all $j=1,\ldots,d$, we conclude that $h^*=0$, which in turn implies that $\{h^*,h_2,\ldots,h_d\}$ is a linearly dependent collection, thereby contradicting our assertion of $\{h^*,h_2,\ldots,h_d\}$ being an orthonormal basis of $\real^d$, which proves our claim.

Having established that there exists $j\in \{1,\ldots,d\}$ with $\langle h^*,f_j\rangle \neq 0$, we encounter two cases. Either $\langle f_j,h^*\rangle > 0$, in which case $HD(tf_j,\pr_n)$ will exhibit an exponential decay as a consequence of $(ii)$ above. Alternatively, $\langle f_j,h^*\rangle < 0$, and then $HD(-tf_j,\pr_n)$ will exhibit an exponential decay, again as a consequence of $(ii)$ above. 

\subsection{Proofs of sample versions}\label{ss:sample-proofs}

All the proofs of statement about sample versions of halfspace depth are presented in this section. Note that all the (in)equalities involving random quantities are to be interpreted as almost sure (in)equalities.

We also note here that the proofs of results related to the sample version of halfspace depths are heavily reliant on the corresponding results for the population version, combined with Theorem~\ref{thm:empHDdecay-vs-popHD}. Specifically, the proofs of Theorems \ref{thm:empirical-hd-decay-tail} and \ref{thm:emp-hdepth-light-decay} are direct consequence of combining Theorem~\ref{thm:empHDdecay-vs-popHD} with Theorem~\ref{thm:hd-decay-tail} and Theorem~\ref{thm:hd-decay-light-tail}$(ii)$, respectively.

\subsubsection{Proofs of Theorem~\ref{thm:empHDdecay-vs-popHD} and Lemma~\ref{prop:hdepth-sup-emp-ratio-gl}}
\label{ss:proofTheoEmpHDdecay}

The proof is articulated in three steps. The first one is based on the observation made in Lemma~\ref{prop:hdepth-sup-emp-ratio-gl}. The second one adapts Lemma~\ref{prop:hdepth-sup-emp-ratio-gl} to the framework given in Theorem~\ref{thm:empHDdecay-vs-popHD}. The third and last step uses Theorem 5.1 in \cite{Alexander1987}.

{\it Step 1.} First notice that Lemma~\ref{prop:hdepth-sup-emp-ratio-gl} follows directly from the observations that: 
\\ If $\displaystyle HD(t_nx,\pr_n) \ge HD(t_nx,\pr)$, then for $\displaystyle\widetilde{H}_{t_nx}$ such that  $\displaystyle\pr(\widetilde{H}_{t_nx}) = HD(t_nx,\pr)$, we can write:
\begin{equation}\label{eq:lemPn>P}
\left|\dfrac{HD(t_n x, \pr_n)}{HD(t_n x, \pr)}-1\right| = \dfrac{HD(t_n x, \pr_n)}{HD(t_n x, \pr)}-1  \le \left(\frac{\pr_n(\widetilde{H}_{t_nx})}{\pr(\widetilde{H}_{t_nx})} - 1 \right) \le \sup_{H \in \calH_{t_nx}} \left|\frac{\pr_n(H)}{\pr(H)} - 1 \right|.
\end{equation}
Similarly, if $HD(t_nx,\pr_n) \le HD(t_nx,\pr)$, then for $H^*_{t_nx}$ s.t. $ \pr_n(H^*_{t_nx})= HD(t_nx,\pr_n)$, we have
$$\left|\dfrac{HD(t_n x, \pr_n)}{HD(t_n x, \pr)}-1\right| =  1- \dfrac{HD(t_n x, \pr_n)}{HD(t_n x, \pr)} \le \left(1 - \frac{\pr_n(H^*_{t_nx})}{\pr(H^*_{t_nx})}  \right) \le \sup_{H \in \calH_{nx}} \left|\frac{\pr_n(H)}{\pr(H)} - 1 \right|.$$

{\it Step 2.} In view of Condition $(C2)$ of Theorem~\ref{thm:empHDdecay-vs-popHD}, we need to consider a subset of $\calH$, hence to adapt Lemma~\ref{prop:hdepth-sup-emp-ratio-gl} to this smaller class, as follows:
\begin{proposition}\label{prop:hdepth-sup-emp-ratio}
    Under Condition $(C2)$ of Theorem~\ref{thm:empHDdecay-vs-popHD}, we have
    \begin{equation}
    \left| \frac{HD(t_nx,\pr_n)}{HD(t_nx,\pr)}-1\right| \le \sup\left\{ \left| \frac{\pr_n(H)}{\pr(H)}-1\right|:\,\,\, H\in \calH_{t_nx},\,\,\pr(H) \ge \gamma_n\right\}.
    \end{equation}
\end{proposition} 
{\bf Proof of Proposition~\ref{prop:hdepth-sup-emp-ratio}.} 

The arguments for the proof of Proposition~\ref{prop:hdepth-sup-emp-ratio} are identical to those used to prove Lemma~\ref{prop:hdepth-sup-emp-ratio-gl}, but taking into account Condition $(C2)$ of Theorem~\ref{thm:empHDdecay-vs-popHD}. Let us consider the first case, when $HD(t_nx,\pr_n) \ge HD(t_nx,\pr)$, for which we have \eqref{eq:lemPn>P}. Then, Condition $(C2)$ implies $HD(t_nx,\pr) = \pr(\widetilde{H}_{t_n x}) \ge \gamma_n$, leading to
$$
\left| \frac{HD(t_nx,\pr_n)}{HD(t_nx,\pr)}-1\right| \le \sup\left\{ \left| \frac{\pr_n(H)}{\pr(H)}-1\right|:\,\,\, H\in \calH_{t_nx},\,\,\pr(H) \ge \gamma_n\right\}.
$$
Next, let us consider the case $HD(t_nx,\pr_n) \le HD(t_nx,\pr)$, where we have
$$\left|\dfrac{HD(t_n x, \pr_n)}{HD(t_n x, \pr)}-1\right| =  1- \dfrac{HD(t_n x, \pr_n)}{HD(t_n x, \pr)} = 1 - \frac{\pr_n(H^*_{t_nx})}{\pr(\widetilde{H}_{t_nx})},$$
where $\widetilde{H}_{t_nx}$ and $H^*_{t_nx}$ are optimal halfspaces for $\pr$ and $\pr_n$, respectively, as defined in the proof of Lemma~\ref{prop:hdepth-sup-emp-ratio-gl} above.
Note that, by definition of halfspace depth,
$$ 
\pr(H^*_{t_nx}) \ge HD(t_nx,\pr) = \pr(\widetilde{H}_{t_nx}).
$$
Now, invoking condition $HD(t_nx,\pr) \ge \gamma_n$, we conclude that $\pr(H^*_{t_nx}) \ge \gamma_n$, which leads to the following upper bound,
$$ \left|\dfrac{HD(t_n x, \pr_n)}{HD(t_n x, \pr)}-1\right| \le 1 - \frac{\pr_n(H^*_{t_nx})}{\pr(H^*_{t_nx})} 
 \le  \sup\left\{ \left| \frac{\pr_n(H)}{\pr(H)}-1\right|:\,\,\, H\in \calH_{t_nx},\,\,\pr(H) \ge \gamma_n\right\}
$$
concluding the proposition. \hfill \qed

{\it Step 3.} This step is based on Theorem 5.1 from \cite{Alexander1987},
which we recall for self containedness in Appendix~\ref{subsec:app2}, to have direct access to the conditions under which it holds.

First, note that the collection $\calH$ of all halfspaces is a VC class, so that Theorem 5.1 from \cite{Alexander1987} holds for the halfspaces. 

Next, since the domain of supremum in Proposition~\ref{prop:hdepth-sup-emp-ratio} is a subset of $\calH$, as a consequence of Alexander's theorem, we can conclude that
\begin{eqnarray}\label{eqn:hdepth-fullsup-reducedsup}
&& \limsup_{n \to \infty} \,\sup \left\{\left|\frac{\pr_n(H)}{\pr(H)}-1\right|: H \in \calH_{t_nx}, \pr(H)\geq \gamma_n\right\}  \nonumber \\
&\le & \limsup_{n \to \infty} \,\sup \left\{\left|\frac{\pr_n(H)}{\pr(H)}-1\right|: H \in \calH, \pr(H)\geq \gamma_n\right\} \nonumber\\
& = & 0\,\,\,\,a.s. 
\end{eqnarray}
Combining this last result \eqref{eqn:hdepth-fullsup-reducedsup} with Proposition~\ref{prop:hdepth-sup-emp-ratio} concludes to Theorem \ref{thm:empHDdecay-vs-popHD}.

\subsubsection{Proof of Theorems~\ref{thm:emp-hdepth-mrv-decay_nodensity} and \ref{thm:emp-hdepth-light-decay}.}
\label{ss:proofTheosTaildecayHD}

Note that the statement of Theorem~\ref{thm:empHDdecay-vs-popHD} also has certain growth conditions on the sequence $t_n$, which in turn are related to $\gamma_n$. \\[1ex]

{\bf Proof of Theorem~\ref{thm:emp-hdepth-mrv-decay_nodensity}.}

{Observe that
\begin{equation*}
    HD(t_n x,\pr_n)=HD(t_n x, \pr)\left[1+\left(\dfrac{HD(t_n x, \pr_n)}{HD(t_n x, \pr)}-1\right) \right].
\end{equation*}
Therefore, we obtain
$$
    \begin{aligned}
     &\sup_{\|x\|=\veps}\left|\frac{HD(t_n x,\pr_n)}{1-\pr(t_n B^d)}- \inf_{\|p\|=1} \int_{\{z:\langle z-x,p\rangle \geq 0\}} \la(z) dz \right|\,\le \\
     & \sup_{\|x\|=\veps}\left\{ \frac{HD(t_n x, \pr)}{1-\pr(t_n B^d)} \left|\dfrac{HD(t_n x, \pr_n)}{HD(t_n x, \pr)}-1 \right|  + \left|\frac{HD(t_n x,\pr)}{1-\pr(t_n B^d)}- \inf_{\|p\|=1} \int_{\{z:\langle z-x,p\rangle \geq 0\}} \!\!\!\!\!\!\!\!\!\!\!\!\!\!\!\la(z) dz \right|\right\}.
    \end{aligned}
$$
    It follows directly from Proposition~\ref{prop:depth_decay_pop_nodensity} that the second term tends to zero.
    Moreover, by Proposition~\ref{prop:depth_decay_pop_nodensity}, for sufficiently large $n$, there exists $c>0$ such that,
$$
HD(t_n x, \pr) > c\, [1-\pr(t_n B^d)].
$$
Combining this last inequality with the given condition $1-\pr(t_n B^d)>\gamma_n$ gives $HD(t_n x, \pr)>c \, \gamma_n$. 

Therefore, it follows from Theorem~\ref{thm:empHDdecay-vs-popHD} that the first term almost surely converges to zero. \hfill }\qed \\[1ex]

{\bf Proof of Theorem~\ref{thm:emp-hdepth-light-decay}.}
The proof follows directly by combining \eqref{eqn:emp-pop-hdepth} with the limit from Theorem~\ref{thm:hd-decay-light-tail}.
\section{Conclusion}
\label{sec:concl}

Much literature has been developed so far on halfspace depth functions, primarily focusing on their properties, such as continuity, convexity, affine equivariance, invariance under (orthogonal) transformations, among others.
Our focus, however, is on the asymptotic properties, particularly exploring its relationship with the tail behaviour of the underlying distribution. 

First, we considered the population side, completing the asymptotics literature to provide a comprehensive understanding and laying the groundwork for the sample side.
We then addressed the same questions regarding the asymptotics when considering the empirical distribution. This is a crucial problem in view of applications, questioning the relevance of these tools when working with samples and remaining within the data hull. In doing so, we complement the work by \cite{Einmahl2015,He2017} in the heavy tail case and that of  \cite{Burr2017} in the light tail case.  

It is worth recalling that halfspace depths do not uniquely identify the underlying probability measure, as recently shown by \cite{Nagy2021a}. However, the characterisation is unique when the measures have finite support, as is the case with empirical measures. This motivated us to further investigate the halfspace depth for samples. 

We derived rates of decay for halfspace depth functions when considering the empirical distribution, in the almost sure sense. 
Additionally, we specified these rates depending on the type of tail behaviour of the measure, whether light or heavy — an important consideration in risk analysis. To further address this practical question, we developed a methodology to identify the nature of the tail behaviour of the underlying measure by analysing the empirical halfspace depth along different directions.

These results are significant both theoretically and practically, as they contribute to a deeper understanding of halfspace depth and help validate its empirical use within different frameworks.

\section*{Acknowledgments}

This work is part of the doctoral thesis of the first author while he was a student at TIFR–CAM, under the supervision of the 2 other authors. It benefited a lot from the mutual visits at ESSEC Business School, Paris, France and at TIFR-CAM, Bangalore, India, respectively. The authors are grateful to the hosting institutions and for funding to the {\it Fondation des Sciences de la Mod\'elisation} (ANR-11-LABX-0023-01), and SERB--MATRICS grant MTR/2020/000629.
   
\bibliographystyle{chicago}
\bibliography{LitSibsankar-202503.bib}

\newpage

\appendix

\begin{center}
\Large{APPENDIX}
\end{center}

\section{MRV framework}\label{app:He2017}

Recall the notion of multivariate regular variation (MRV). A random vector $\mathbf{X}$ is multivariate regularly varying if there exists
 a probability measure $\nu$ on the unit hypersphere and a (tail) index $\alpha>0$ such that
$$
\lim_{t\to \infty}\frac{\Pr\left(||\mathbf{X}||\ge tx,\,\mathbf{X}/||\mathbf{X}||\in B\right)}{\Pr\left(||\mathbf{X}||\ge t\right)}=x^{-\alpha}\nu(B),
$$
for every $x>0$ and Borel set $B$ in the unit hypersphere such that $\nu(\partial B)=0$.

An equivalent statement at the density level corresponds to the two convergence conditions given in Corollary~\ref{cor:theo-emp-hdepth-mrv}.

Let us recall Proposition~2 in \cite{He2017}, using our notations.

{\begin{proposition}[\textnormal{Proposition~2 in \cite{He2017}}]
\label{prop:depth_decay_pop_nodensity}
    Let $\pr$ be a probability measure defined on $(\real^d,\calB(\real^d))$ that satisfies all the conditions of Theorem \ref{thm:emp-hdepth-mrv-decay_nodensity}. 
    Then, for any $\veps > 0$, we have 
    \begin{equation}
        \lim_{t \to \infty} \sup_{\|x\| \geq \veps}\left|\frac{HD(tx,\pr)}{1-\pr(tB^d)} - HD(x, \nu)\right|=0.
    \end{equation}      
\end{proposition}
}

Assuming that $\pr$ has a density $f$, we can adapt Theorem \ref{thm:emp-hdepth-mrv-decay_nodensity} as follows: 
\begin{corollary}
\label{cor:theo-emp-hdepth-mrv}
    Let $\pr$ be a probability measure on $(\real^d,\calB(\real^d))$ with density $f$ continuous on a neighborhood of infinity, such that the map $y\mapsto \|y\|^d f(y)$ is bounded in every compact neighbourhood of the origin, and there exist $\la:\real^d \to \real^{+}$ and $V\in RV_{-\alpha}$, with $\alpha>0$, such that
    \begin{equation*}
        \left|\frac{f(ty)}{t^{-d}V(t)} -\la(y)\right| \underset{t\to\infty}{\longrightarrow} 0,\,\,\,\,\,\forall y\neq 0,
        \quad \text{and} \quad
        \sup_{\|y\|= 1}\left| \frac{f(ty)}{t^{-d}V(t)} -\la(y)\right| \underset{t\to\infty}{\longrightarrow} 0.
    \end{equation*}
Let $\{X_n\}_{n\ge 1}$ be an i.i.d sample drawn from $\pr$. Then, we have,
    \begin{equation}
        \lim_{n \to \infty} \sup_{{\|x\|=1} }\left|\frac{HD(t_n x,\pr_n)}{V(t_n)} - \inf_{\|p\|=1} \int_{\{z:\langle z-x,p\rangle \geq 0\}} \la(z) dz \right| = 0, \,\,\,\,\,\,\,\text{a.s.}
    \end{equation}
whenever $(t_n)\nearrow \infty$ such that $V(t_n)>\, \gamma_n$ (for any large $n$), and $(\gamma_n)$ satisfies Conditions $(C1a)$ and $(C1b)$ given in Theorem \ref{thm:empHDdecay-vs-popHD}.
\end{corollary}
%
\section{Alexander's framework}\label{subsec:app2}
 
For the paper to be self-contained, we shall recall some concepts that are fundamental to understanding convergence of empirical processes, as used in establishing convergence of empirical halfspace depth.
\vskip 0.2cm

{\bf Vapnik--\u{C}ervonenkis (VC) class}

For simplicity, we shall consider the measure space $(\real^d,\calB(\real^d))$. A class of subsets $\calC$ of $\real^d$ is said to be of VC class, if 
$$\sup\{\text{card}\{F\cap C: C\in \calC\}: \text{card}(F)=n,\, F\subset \real^d\} < 2^n$$
for some $n\ge 1$.

Note that the set $\calH$ of all halfspaces of $\real^d$ is a VC class (see \cite{Alexander1987} and references therein for more details).

{\bf Capacity function}

For the VC class of sets $\calC$ and a measure $\pr$ on $(\real^d,\calB(\real^d))$, we define
\begin{eqnarray*}
\calC_t & = & \left\{C\in\calC: \pr(C)(1-\pr(C))\le t,\,\,\pr(C)\le \frac12\right\} \,\cup \\
&& \qquad\qquad\left\{C^c\in\calC: \pr(C)(1-\pr(C))\le t,\,\,\pr(C)> \frac12\right\}; \\
E_t & = & \bigcup_{C\in \calC_t} C\,; \qquad
a_t  \,=\,  \pr(E_t) \vee t\,; \qquad
g_c(t) \, = \, \frac{a(t)}{t} \,.
\end{eqnarray*}
The function $g_c$ is called the capacity function corresponding to $\calC$ and $\pr$. \cite{Alexander1987} discusses and analyses various properties of $g_c$. Quoting Alexander, `{\it $g_c$ can be thought of roughly as the number of disjoint sets of size $t$ which will `fit' in $\calC$: $a(t)$ is the space available, and $t$ is the approximate space needed for each set $C$ with $\pr(C)\simeq t$.}'

Note that, by definition of $g_c$, we always have 
$$
g_c(t)\le 1/t, \,\, \text{for any } t.  
$$
Specifically, in {\bf Example 3.6 of \cite{Alexander1987}}, the author shows that, in the Gaussian case, 
$$
 g_c(t) \sim K \left( \log (1/t)\right)^{(d-1)/2}, \quad \text{as} \quad t\to 0, \quad (K \;\text{being some constant}),
 $$ 
 if $\calC$ is taken to be the set of all halfspaces in $\real^d$ and $\pr$ is set to be the $d$-dimensional standard Gaussian. 

 Now, turning to MRV setting, we can notice that 
$$
g_c(t)=O(1),
$$ 
a result already used in \cite{Einmahl2015,He2017}. Let us briefly prove this, as it was not explicitly shown in the quoted papers.

When $\mathcal{C} = \mathcal{H}$, the set $\mathcal{C}_t$ can be expressed for sufficiently small $t$ as  
$$
\mathcal{C}_t = \left\{C \in \mathcal{C} : \pr(C)(1-\pr(C)) \leq t \right\}.
$$  
Since $\pr(C)(1-\pr(C)) \sim \pr(C)$ for small values of $\pr(C)$, it follows that  
$$
\mathcal{C}_t = \left\{C \in \mathcal{C} : \pr(C) \leq t \right\}.
$$  
Moroever, using Equation (31) of \cite{Einmahl2015}, we know that for MRV distributions, $\pr(E_t) = O(t)$ as $t \to 0$. 
Therefore, by the definition of $g_c(t)$, it follows that $g_c(t) = O(1)$ as $t \to 0$.
\\[1ex]

Finally, let us state Alexander’s result, which sets the stage for our results on the asymptotics of halfspace depth.

{\bf Theorem 5.1 in \cite{Alexander1987}}

    Let $\pr$ be a probability measure defined on $(\real^d,\calB(\real^d))$, and let $g_c$ be capacity function corresponding to $\pr$ as defined in \cite{Alexander1987}. Consider a sequence $\gamma_n$ satisfying the following conditions:
    $$    
    n^{-1}\log (g_c(\gamma_n)) = {\textrm o}(\gamma_n) \qquad \text{and}\qquad n^{-1}\log \log n = {\textrm o}(\gamma_n).
    $$
    Then,
   $$
    \limsup_{n \to \infty} \,\,\,\sup \left\{\left|\frac{\pr_n(C)}{\pr(C)}-1\right|: C \in \calS, \pr(C)\geq \gamma_n\right\}=0\,\,\,\,a.s.
  $$

\section{Illustration of the theoretical results and discussion} 
\label{app:illustration}

To provide a clearer understanding of the theoretical findings, we turn to empirical illustrations. These illustrations explore both light and heavy-tailed distributions.

Halfspace depths are computed using the algorithms proposed in \cite{Dyckerhoff2021} and implemented in the data-depth Python library.


First, we illustrate  Theorems~\ref{thm:emp-hdepth-mrv-decay_nodensity}~\&~\ref{thm:emp-hdepth-light-decay}, considering the asymptotic behaviour of halfspace depth (or Tukey depths) for light and heavy tailed distributions. 
As examples, we consider Gaussian and Pareto distributions for comparison, as we have characterised the rate of convergence according to the tail behaviour.
\begin{figure}[h]
    \centering
    \begin{minipage}{0.48\textwidth}
        \includegraphics[width=0.9\linewidth]{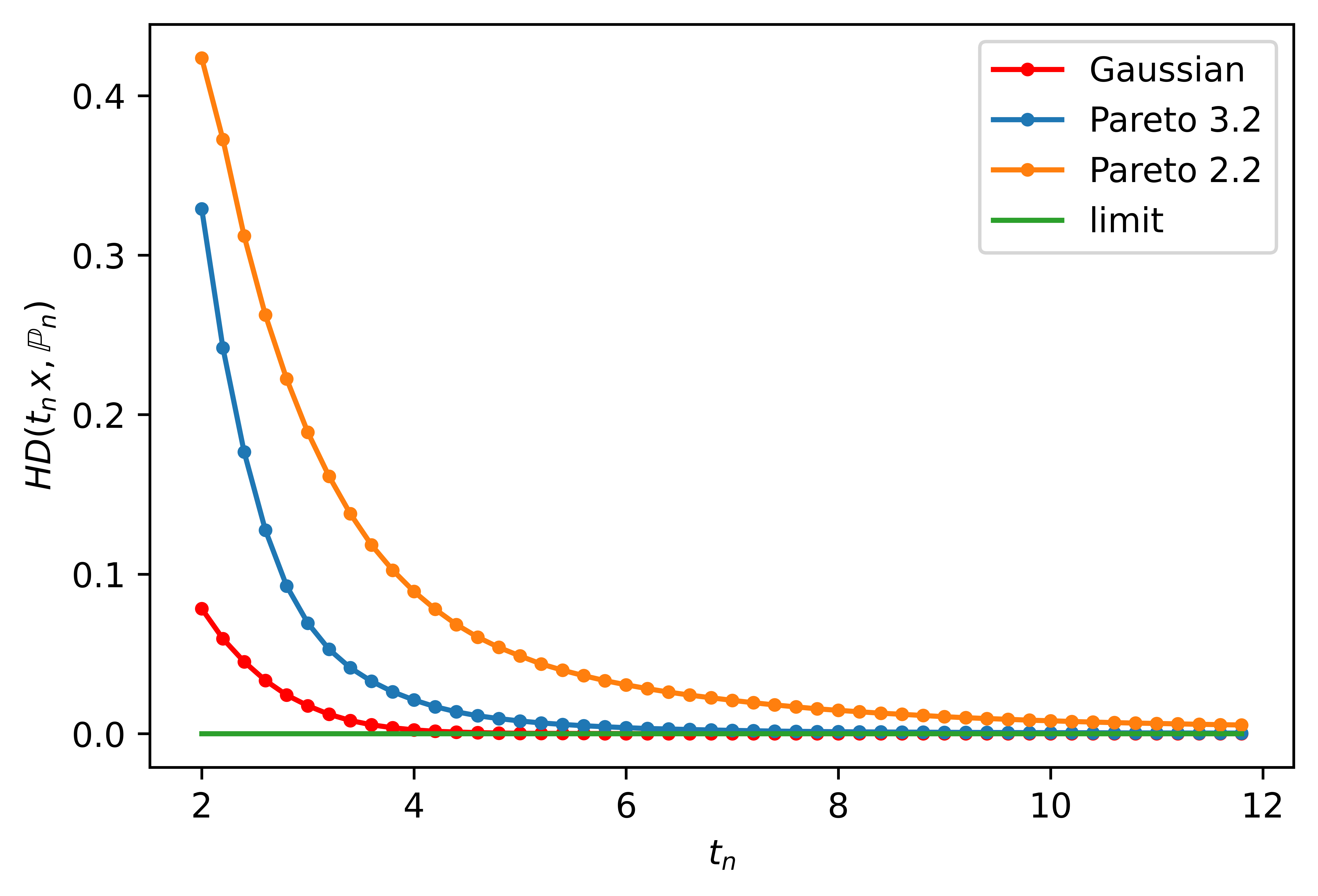}
    \end{minipage}
    \hfill 
    \begin{minipage}{0.48\textwidth}
       \includegraphics[width=0.9\linewidth]{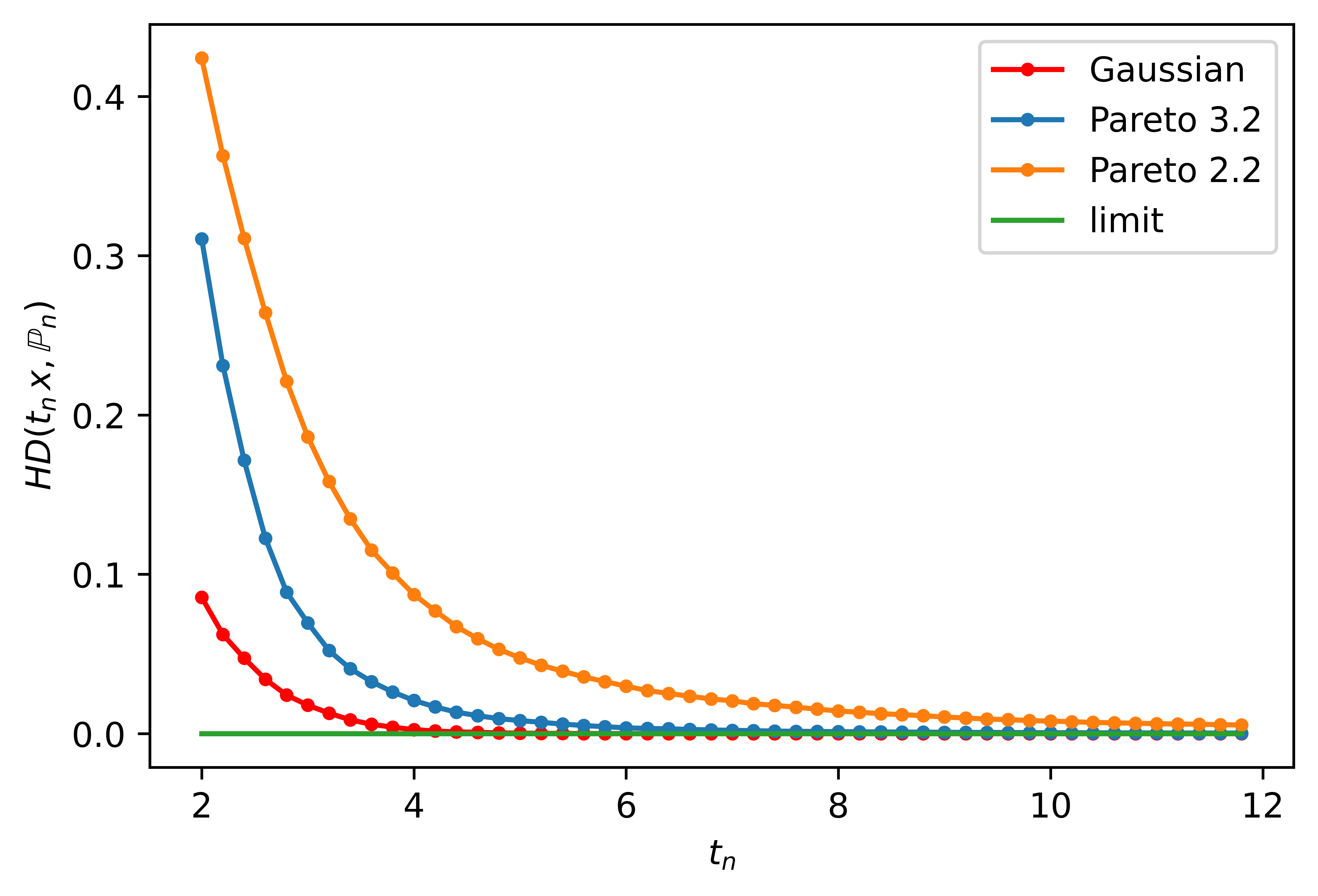}
    \end{minipage}   
\parbox{330pt}{\caption{\sf \small Tukey depths are computed at points in direction $x=(1,1)$ and given in terms of $(t_n)$ growing linearly in $n$ ($t_n=1.8+n.10^{-4}$, with $n=10^5 k/50$, $k=1,2\cdots,50$). Samples are taken from independent bivariate Pareto with parameter $2.2$ and $3.2$, respectively, and Gaussian distribution with diagonal covariance matrix $diag(2,2)$. Number of observations is $10^5$. Left plot: fixed sample. Right plot: growing sample, with a partition of 50}
\label{fig:ratecvdepth-asymp-linearTime}}
\end{figure}

In Figure~\ref{fig:ratecvdepth-asymp-linearTime}, we plot $HD(t_n\,x, \pr_n)$ as a function of $t_n$, 
choosing for $x$ the direction $(1,1)$, $\pr_n$ coming from, respectively, bivariate standard Gaussian and Pareto($\delta$) distributions, with independent components, and  $\delta=1.9$, $2.2$ and $3.2$, to span the spectrum from very heavy to moderately heavy tail. 
We compare $t_n$ when taking a fixed sample (left plot) and a linearly growing one (right plot).
Comparing the different depths according to the type of distributions, from very light to moderate heavy (with second moment but no third one), we clearly observe a different rate of convergence towards $0$. The heavier is the distribution, the slower is the convergence. 

\begin{figure}[h]
    \centering
    \begin{minipage}{0.325\textwidth}
        \includegraphics[width=0.9\linewidth]{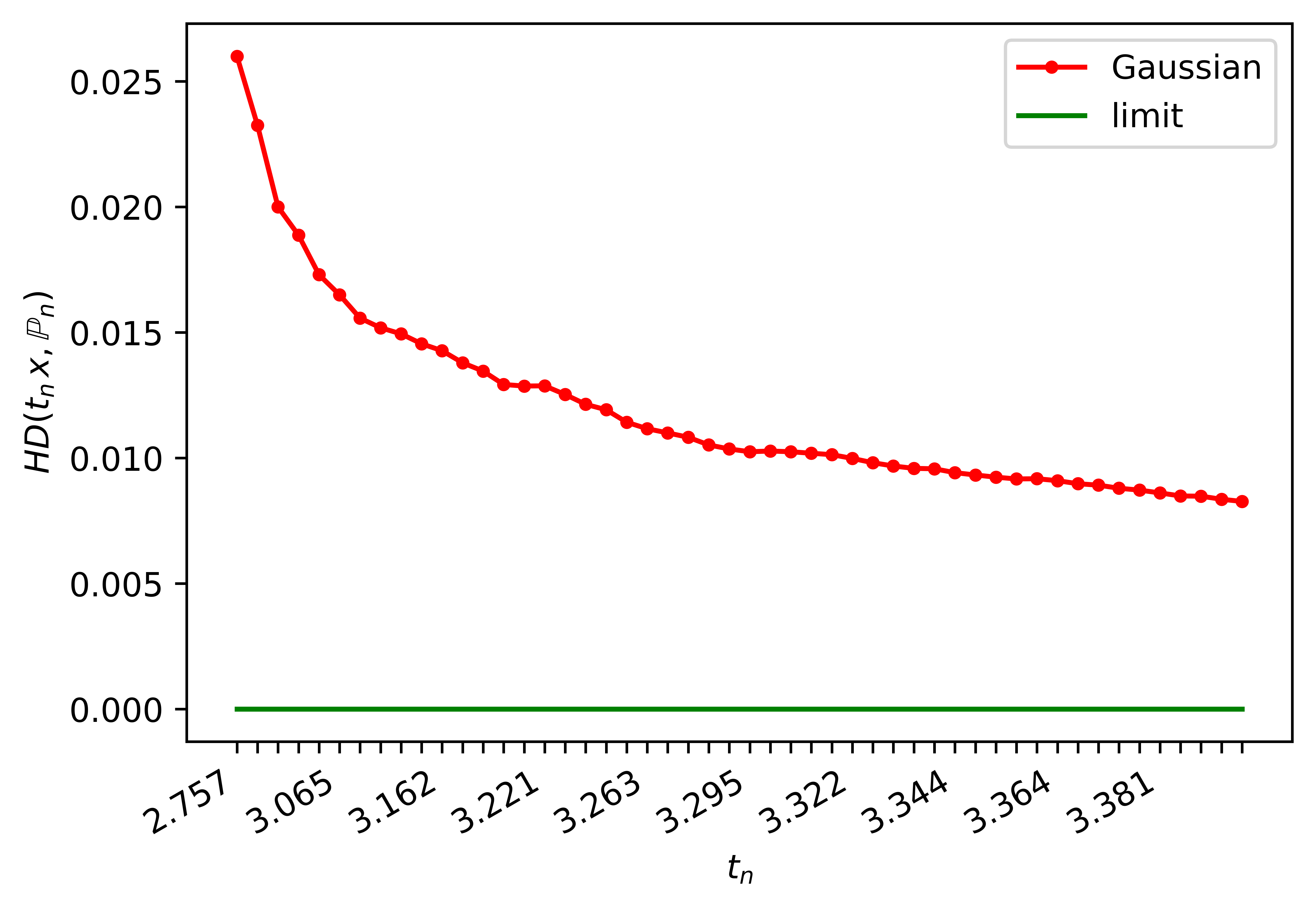}
    \end{minipage}
    \hfill 
    \begin{minipage}{0.325\textwidth}
        \includegraphics[width=0.8\linewidth]{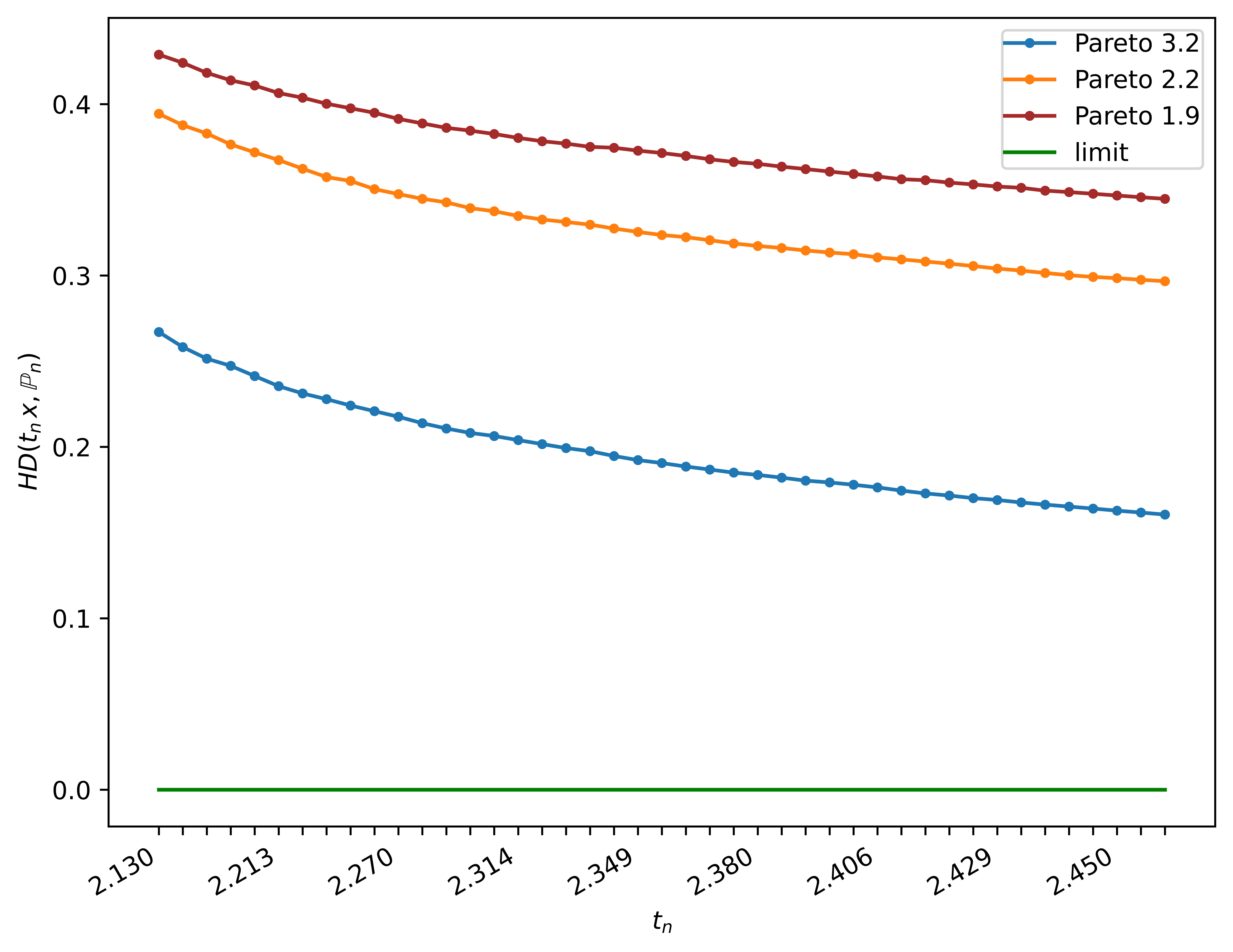}
    \end{minipage}
    \hfill
    \begin{minipage}{0.325\textwidth}
        \includegraphics[width=0.9\linewidth]{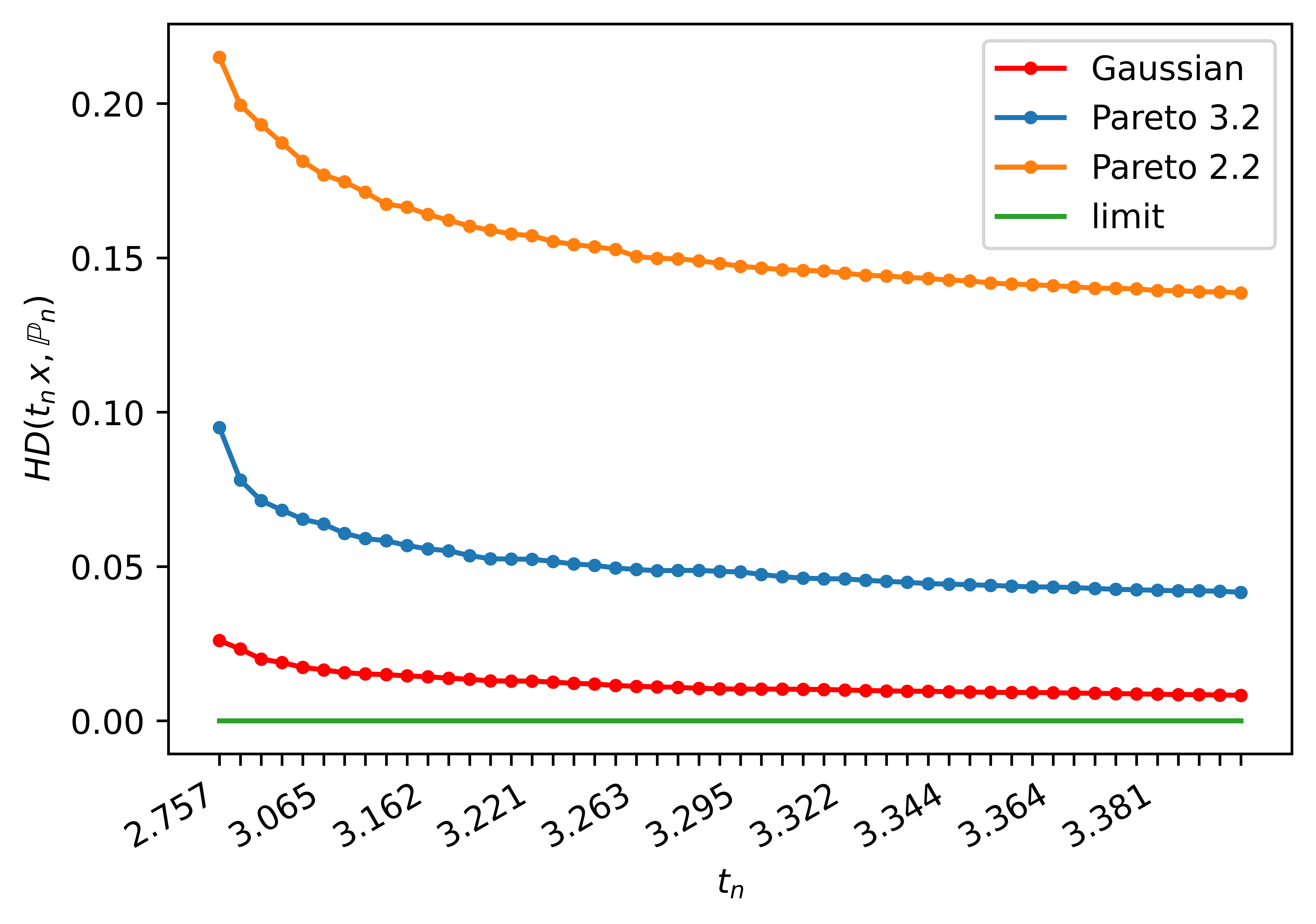}
    \end{minipage}    
\parbox{340pt}{\caption{\sf \small Tukey depths are computed at points in direction $x=(1,1)$ given in terms of $(t_n)$. Samples are taken from independent bivariate Pareto with parameter $1.9$, $2.2$ and $3.2$, respectively, and Gaussian distribution with diagonal covariance matrix $diag(2,2)$. Number of observations is $10^5$. Left plot: Halfspace depth for the Gaussian sample. Middle plot: Comparing the halfspace depth behaviours when considering Pareto($\delta$) samples, varying $\delta$, $t_n$ corresponding to $t_n(\text{Pareto}(3.2))$. Right plot: Halfspace depth behaviours for Gaussian and Pareto($\delta$) samples (choosing $\delta>2$) and for $t_n=\sqrt{\log n}$ (Gaussian).}
\label{fig:ratecvdepth-asymp}}
\end{figure}

Next, building on the rate of convergence found in the light tail (see Theorem~\ref{thm:emp-hdepth-light-decay} and Example~\ref{exple:exp-normal}) and the heavy tail case (see Theorem~\ref{thm:emp-hdepth-mrv-decay_nodensity}), we plot $HD(t_n x,\pr_n)$ as a function of $t_n$. Given the very different speeds of convergence obtained for the light versus heavy tails, we first give a plot for the Gaussian sample only, then a plot for Pareto($\delta$) samples with varying $\delta$, so that we can appreciate the different behaviour and convergence depending on the heaviness. The sequence $\{t_n\}$ is chosen according to the type of distribution. For the Gaussian case (see Figure~\ref{fig:ratecvdepth-asymp}, left plot), $t_n=\sqrt{\log n}$ (choosing $\beta=1/2$ in Example~\ref{exple:exp-normal},(b)). For Pareto distributions (see Figure~\ref{fig:ratecvdepth-asymp}, middle plot), we consider 
$t_n=n^{-\frac{\beta}{2\delta}}$, with $0<\beta<1$ also chosen as $1/2$ and the Pareto parameter $\delta$ corresponding to the less heavy, {\it i.e.} $\delta=3.2$ (since the lighter the tail, the faster the convergence towards $0$). Finally, we provide a last plot (see Figure~\ref{fig:ratecvdepth-asymp}, right plot) comparing the Gaussian and Pareto cases, choosing the scaling $t_n$ associated with the Gaussian distribution, for better visualizing the difference of behaviours and speeds of convergence.
The three plots given in Figure~\ref{fig:ratecvdepth-asymp} highlight the difference of rates of decay of the halfspace depths according to the tail behaviour of the measure. The left and middle plots point out the fast convergence of halfspace depth for the Gaussian sample (decreasing from 2.6\% to less than 1\% (0.83\%) on the given range for $t_n(\text{Gaussian})$), and the impact of the heaviness for the Pareto samples, with a decrease from $26.7\%$ to $16\%$ on the given range for $t_n(\text{Pareto}(3.2))$ for the Pareto with 3rd moment, from 39\% to 30\% for the Pareto(2.5), while from 43\% to 35\% for the heaviest Pareto (with no 2nd moment), hence a very slow decrease compared with Pareto(3.2). The third plot allows for a direct comparison between light and heavy tails, considering the Gaussian scaling for $t_n$; the relation between the rate of decay of the halfspace depth to $0$ and the tail behaviour becomes even more obvious. 

Note that it would have been nice to look at the convergence towards $1$ of the normalized halfspace depth function $HD(t_n x,\pr_n)/N(t_n)$ (rather than $HD(t_n x,\pr_n)$), as a function of $t_n$ , where $N(t_n)$ corresponds to the speed of convergence, namely of order $N(t_n)=n^{-\eta}$ for the Gaussian case (see Example~\ref{exple:exp-normal}(b)) and $N(t_n)=1-\pr(t_n B^d)$ as defined in Theorem~\ref{thm:emp-hdepth-mrv-decay_nodensity}  for the Pareto one. Nevertheless, to observe something informative in terms of convergence, it would require a large number of observations (more than $10^{20}$), which is computationally not feasible with the R-package we are using. We conjecture that a possible way to circumvent this computational hurdle would be to use the geometry of isoquantile (isodepth) contours, as they secrete immense amount of information about the underlying distribution. 

\section{Annexure to the real data example} 
\label{app:OLR}

\begin{figure}[h]
    \centering
    \begin{minipage}{0.44\textwidth}
        \includegraphics[width=0.9\linewidth]{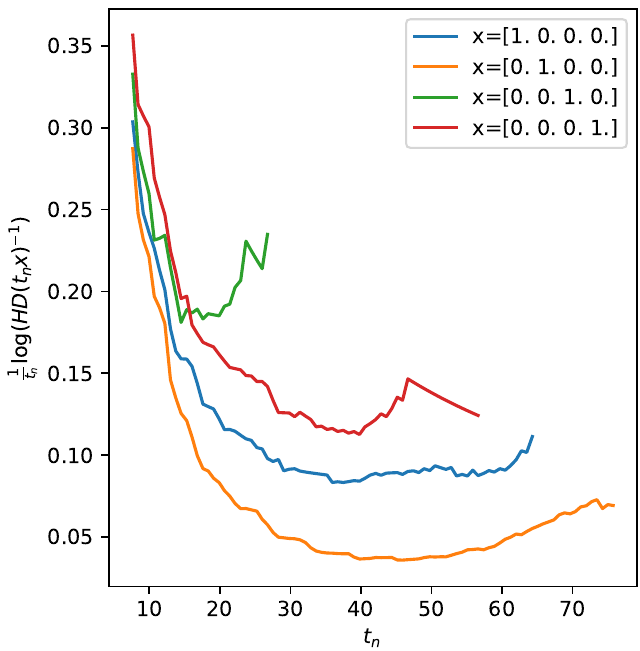}
    \end{minipage}
    \hfill
    \begin{minipage}{0.44\textwidth}
       \includegraphics[width=0.9\linewidth]{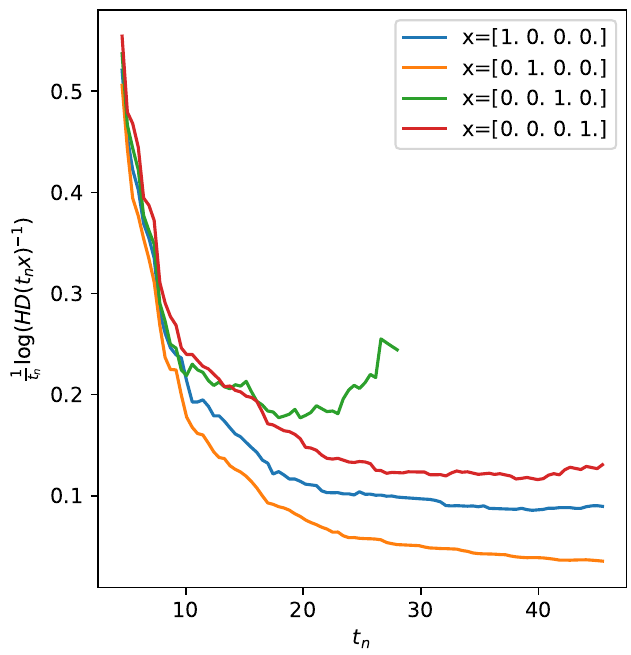}
    \end{minipage}   
\parbox{330pt}{\caption{\sf \small The left and right plots correspond to $\displaystyle y_n := \log\left(\left( {HD(t_n x, \pr_n)} \right)^{-1}\right)/{t_n} $ plotted for the OLR data in $e_1, e_2, e_3$ and $e_4$ directions against $t_n$, with $t_n = 5n/100$ for the left plot, and $t_n=3n/100$ for the right plot, for $n=\lfloor 17.08 k\rfloor$, for $k=10,\ldots,100$.}
\label{fig:realData-HDrate-RightTail}}
\end{figure}

Here, we plot $HD(t_n x,\pr_n)$ for $x=e_1,e_2,e_3,e_4$ in Figure \ref{fig:realData-HDrate-RightTail}. It is clearly seen from the plots that Prayagraj Winter has the lightest tail, and Jaisalmer Summer the heaviest tail among the four. It can be interpreted as the OLR at Prayagraj in winter is very rarely higher than its seasonal mean, whereas the OLR at Jaisalmer in summer has slightly higher chance of taking values higher than its seasonal mean.

\end{document}